\documentclass[12pt,twoside]{amsart}
\usepackage{amssymb}

\nonstopmode

\numberwithin{equation}{section}
\font\tengothic=eufm10 scaled\magstep 1
\font\sevengothic=eufm7 scaled\magstep 1
\newfam\gothicfam
      \textfont\gothicfam=\tengothic
      \scriptfont\gothicfam=\sevengothic
\def\goth#1{{\fam\gothicfam #1}}

\newtheorem{theorem}{Theorem}[section]
\newtheorem{lemma}[theorem]{Lemma}
\newtheorem{proposition}[theorem]{Proposition}
\newtheorem{corollary}[theorem]{Corollary}

\newtheorem{question}[theorem]{Question}

\theoremstyle{definition}
\newtheorem{definition}[theorem]{Definition} 
\newtheorem{remark}[theorem]{Remark}
\newtheorem{example}[theorem]{Example}

\newtheorem{notation}[theorem]{Notation}

\newcommand{\codim}{\operatorname{codim}}

\newcommand{\Hom}{\operatorname{Hom}}
\newcommand{\Ext}{\operatorname{Ext}}
\newcommand{\reg}{\operatorname{reg}}

\newcommand{\rank}{\operatorname{rank}}

\newcommand{\depth}{\operatorname{depth}}

\newcommand{\im}{\operatorname{im}}

\newcommand{\TT}{\operatorname{Tor}}
\newcommand{\tor}{\operatorname{tor}}

\newcommand{\Proj}{\operatorname{Proj}}
\newcommand{\proj}[1]
{ \mathchoice
           { {\mathbb P}^{#1} }
           { {\mathbb P}^{#1} }
           { {\mathbb P}^{#1} }
           { {\mathbb P}^{#1} }
         }

\newcommand{\ffi}{\varphi}

\newcommand{\acm}{arithmetically Cohen-Macaulay }
\newcommand{\aG}{arithmetically Gorenstein }
\newcommand{\WLP}{weak Lefschetz property }

\newcommand {\PP}{\mathbb{P}}

\begin{document}
\title[  Gorenstein schemes and simplicial polytopes]{Reduced
  arithmetically Gorenstein schemes and 
simplicial polytopes with maximal Betti numbers}

\author[J.\ Migliore, U.\ Nagel]{J.\ Migliore$^*$ 
\and 
U.\ Nagel$^{**}$ }

\address{Department of Mathematics,
        University of Notre Dame,
        Notre Dame, IN 46556,
        USA} 

\email{Juan.C.Migliore.1@nd.edu} 

\address{Fachbereich Mathematik und Informatik,
  Universit\"at-Gesamthochschule 
Paderborn, D--33095 Paderborn, Germany} 

\email{uwen@uni-paderborn.de} 

\thanks{$^*$ Partially supported  by the Department of Mathematics
of the  University of Paderborn\\
$^{**}$ Partially supported  by the Department of Mathematics of the
  University of Notre Dame \\
2000 {\em Mathematics Subject Classification:} Primary 13C40, 13D02, 13D40, 
14M05, 14M06, 14N20; \\ Secondary 13H10, 52B12, 52B11 }


\begin{abstract}
An SI-sequence is a finite sequence of positive integers which is
symmetric, unimodal and satisfies a certain growth condition.  These are
known to correspond precisely to the possible Hilbert functions of
Artinian Gorenstein algebras with the Weak Lefschetz Property, a property
shared by most Artinian Gorenstein algebras.   Starting
with an arbitrary SI-sequence, we construct a reduced, arithmetically
Gorenstein configuration $G$ of linear varieties of arbitrary dimension
whose Artinian reduction has the given SI-sequence as Hilbert function
and has the Weak Lefschetz Property.  Furthermore, we show that $G$ has
maximal graded Betti numbers among all arithmetically Gorenstein
subschemes of projective space whose Artinian reduction has the Weak
Lefschetz Property and the given Hilbert function. As an application we
show that over a field of characteristic zero every set of simplicial
polytopes with fixed $h$-vector contains a polytope with maximal graded
Betti numbers.
\end{abstract}


\maketitle

\tableofcontents


 \section{Introduction} \label{intro}

This paper addresses three fundamental questions about arithmetically
Gorenstein subschemes of projective space.  First, we consider the question of
possible Hilbert functions that can occur (in arbitrary codimension).  Second,
we consider the possible graded Betti numbers that can occur for the minimal
free resolution of the homogeneous ideal of such a subscheme.  In particular,
we are interested in the problem of whether, among the arithmetically
Gorenstein schemes with fixed Hilbert function, there is one with maximal
graded Betti numbers.   Most importantly, we are interested in the question of
liftability from the Artinian case: which properties of Artinian Gorenstein
ideals lift to properties of the ideal of a {\em reduced}, arithmetically
Gorenstein subscheme of projective space?  (Of course without the ``reduced''
requirement all properties can be lifted by considering cones.)  In
particular, which Hilbert functions and which graded Betti numbers  that occur
at the Artinian level also occur for reduced arithmetically Gorenstein
schemes? We also consider these questions for the even more special reduced,
monomial, Gorenstein ideals occurring as Stanley-Reisner ideals of
simplicial polytopes. 

We will give complete answers to these questions for a large subset
of the set of all arithmetically Gorenstein subschemes, namely for the
ones whose Artinian reductions have the so-called {\em weak Lefschetz
property}.   The weak Lefschetz property is a very important notion for
Artinian Gorenstein algebras; it says that the homomorphism induced between
any two consecutive components by multiplication by a general linear form has
maximal rank.   Watanabe \cite{watanabe} has shown that ``most'' artinian
Gorenstein algebras have this property. {\em From now on we say that an \aG
scheme $X$ has the \WLP if there is an Artinian reduction of  $X$ having the
weak Lefschetz property.}

Codimension three Gorenstein ideals are quite well understood,
thanks primarily to the structure theorem of Buchsbaum and Eisenbud
\cite{BE}.  Using this result, Diesel \cite{diesel} gave a description of all
possible sets of graded Betti numbers (hence all possible Hilbert functions)
for {\em Artinian} Gorenstein ideals, leaving open the question of which of
these could occur for reduced, arithmetically Gorenstein subschemes of
projective space.  Geramita and the first author solved this problem by
showing that any set of graded Betti numbers allowed by Diesel in fact occurs
for a reduced set of points in $\proj{3}$, a {\em stick figure} curve in
$\proj{4}$, and more generally a ``generalized stick figure'' in
$\proj{n}$.    As a consequence, they showed that any arithmetically
Gorenstein subscheme of $\proj{n}$ ($n \geq 3$) specializes to a good linear
configuration with the same graded Betti numbers.

What can be said in higher codimension?  In codimension $\geq 4$ it is not
even known precisely which Hilbert functions arise for Artinian Gorenstein
ideals, and there is no analog to Diesel's work for graded Betti numbers.
Still, for Hilbert functions there are some results in \cite{boij}, \cite{BL}
and \cite{harima}, for example, and for resolutions there are important results
in \cite{GHS1} which we will describe below.  These results form the starting
point for our work.

Although a complete classification of possible Hilbert functions for Artinian
Gorenstein ideals of codimension $\geq 4$ is not known, a very large class of
Artinian Gorenstein Hilbert functions have been constructed by Harima
\cite{harima}, namely those which are so-called SI-sequences.  Roughly, this
property says that the ``first half'' of the Hilbert function (which is a
finite sequence of integers) is a differentiable O-sequence, in the sense of
Macaulay \cite{mr.macaulay} and Stanley \cite{stanley}.  It is known that not
all Hilbert functions of Artinian Gorenstein ideals in codimension $\geq 5$
have this property, but it is an open question whether they all do in
codimension 4.  Certainly they all do in codimension 3.  Harima's approach was
via liaison and sums of linked ideals, to which we will return shortly.

Prior to Harima's work, it was already shown by Billera and Lee \cite{BL2} and
Stanley \cite{stanley2} that an $h$-vector is the $h$-vector of a simplicial
polytope if and only if it is an SI-sequence.  As observed by Boij
\cite{boij}, this implies that for any given SI-sequence there exists a
reduced Gorenstein algebra whose
$h$-vector is that sequence.  However, as Harima points out, Stanley's
methods used 
in \cite{stanley2} involve hard results about toric varieties and topology
and apply in characteristic zero only. Thus it is 
worth giving a different proof.  The
first part of this paper gives a new proof of this fact, 
providing a ``lifting'' of Harima's Artinian result (but not his proof).  
Our approach is similar to his, but with some important differences which
will be described shortly.  More precisely, our first main result is the
following:

\begin{theorem} \label{1stmain}
Let $\underline{h}= (1,c,h_2, \dots h_{s-2},c,1)$ be an SI-sequence.
Then, for every integer 
$d \geq 0$, there is a reduced arithmetically Gorenstein scheme $G\subset
\PP_K^{c+d}$ of dimension $d$, with the weak Lefschetz property, whose
$h$-vector is $\underline h$, 
provided the field $K$ contains sufficiently many elements. 
\end{theorem}

Note that it is easy to see that the $h$-vector of an \aG with the \WLP is
an SI-sequence. Thus combining with Watanabe's result, Theorem \ref{1stmain}
provides a characterization of the Hilbert functions of ``most" reduced \aG
schemes (especially if Question \ref{reduced-WLP} below has a positive answer).

  We were not able to mimic the approach of Harima
to prove Theorem \ref{1stmain} (see Remark \ref{compare-harima-1} and Remark
\ref{compare-harima-2}).  Such an approach would amount to adding the ideals of
certain reduced \acm subschemes of projective space which are linked {\em by a
complete intersection}, thus forming a reducible arithmetically Gorenstein
subscheme of codimension one more with a ``bigger'' Hilbert function than the
desired one, and then removing components to obtain an arithmetically
Gorenstein subscheme with the desired Hilbert function.
Instead, our approach is to link using arithmetically Gorenstein ideals
(i.e.\ G-links, in the terminology of \cite{KMMNP}) rather than complete
intersections (i.e.\ CI-links).   We refer to \cite{migliore} for
background on liaison.  To the best of our knowledge this is one of the first
occurrences in the literature of using Gorenstein liaison to
construct interesting subschemes of projective space.  In fact we show in
Remark \ref{ci-not-enough} that complete intersections do not suffice to
obtain all the Hilbert functions.

The main features of our construction are the
following.

\begin{itemize}
\item A simple calculation shows that the desired Hilbert function can be
obtained (at least numerically) by linking an \acm scheme with the ``obvious''
Hilbert function (namely the one suggested by the ``first half'' of the desired
$h$-vector) inside an arithmetically Gorenstein scheme with ``maximal'' Hilbert
function.  This calculation is given in section \ref{tie-together}.  The hard
part is to verify that suitable schemes can actually be found to carry this
out.

\item We observe that if a Gorenstein scheme is a generalized stick figure
then any link which it provides will be geometric, and the sum of linked
ideals will be reduced.

\item We give a useful construction of \acm schemes:  Given a nested sequence
$V_1 \subset V_2 \subset \cdots \subset V_r$ of generically Gorenstein \acm
subschemes of the same dimension, the union of generally chosen hypersurface
sections of the $V_i$ is again an  \acm subscheme of $\proj{n}$, and its
Hilbert function is readily computed.  This construction is central to the
argument, and is of independent interest (cf.\ for instance \cite{MN4}).

\item Using the  construction just mentioned, we show in section
\ref{max-hf-sect} how to make \acm generalized stick figures with ``maximal''
Hilbert function in any codimension.  We give a precise primary decomposition
of these configurations.  We also show how a generalized stick figure
with arbitrary Hilbert function can be found as a subconfiguration of such a
``maximal'' one.

\item  We show how to construct arithmetically Gorenstein generalized stick
figures with ``maximal'' Hilbert function  which contain the  ``maximal''
generalized stick figures just mentioned.  Note that the Artinian reductions
of our \acm schemes with maximal Hilbert functions are compressed
Cohen-Macaulay  algebras in
the sense of Iarrobino \cite{iarrobino}, \cite{IK}.  We give a precise primary
decomposition of these configurations.

\item Tying things together, we show that this set-up allows us to carry out
the program described in the first step above.

\end{itemize}

As a result,  we have a purely combinatorial way of getting
arithmetically Gorenstein schemes with the desired Hilbert functions.  As a
consequence, we have the interesting fact that our links are obtained by
beginning with the ``linking'' arithmetically Gorenstein scheme and then
finding {\em inside it} the subschemes which are to be linked.  This is the
opposite of the usual application of liaison.

A more precise statement of Harima's result is that {\em a sequence is an
SI-sequence if and only if it is the Hilbert function of an Artinian
Gorenstein algebra with the weak Lefschetz property}.
In order to check that the arithmetically Gorenstein schemes that we construct
have the \WLP we need a higher-dimensional analogue of the \WLP which does not
refer to an Artinian reduction. To this end we introduce in Section
\ref{subspace-sect} the {\em subspace property}. This property is even defined
for \acm schemes. We show that this property is preserved under hyperplane
sections (including Artinian reductions) and that, in the case of an Artinian
Gorenstein algebra, it is equivalent to the weak Lefschetz property. Moreover,
we prove that the \aG schemes constructed above have the subspace property,
finishing the argument for the proof of Theorem \ref{1stmain}.

Furthermore, the subspace property allows us to address a question posed to us
by Tony Geramita: Harima showed that every SI-sequence is the Hilbert function
of some Artinian Gorenstein algebra with the weak Lefschetz property, and
Watanabe \cite{watanabe} showed that ``most'' Artinian Gorenstein algebras have
this property.  However, the question remained whether
reduced arithmetically Gorenstein subschemes of projective space possess the
weak Lefschetz property.  Since we do not know which Artinian Gorenstein
algebras lift to reduced subschemes, it was just possible that these Artinian
reductions are special enough that they avoid the open set of algebras with
the weak Lefschetz property.  However, it follows immediately from the facts
above that this is not the case.  This also follows from work of Stanley
\cite{stanley2} in case the ground field has characteristic zero.

The last part of the paper addresses the question of resolutions.  An
important result in this direction was obtained in \cite{GHS1}.  There the
authors showed that for many, but not all, SI-sequences $\underline{h}$ they
can produce an Artinian Gorenstein algebra $A$ with Hilbert function
$\underline{h}$ possessing the weak Lefschetz property, and such that among
Artinian Gorenstein algebras with the weak Lefschetz property and Hilbert
function $\underline{h}$, $A$ has maximal graded Betti numbers.  Our second
main result of this paper is a generalization of this theorem for reduced
arithmetically Gorenstein schemes of any dimension (since the Artinian
reduction has the weak Lefschetz property and {\em every} SI-sequence can be
produced):

\begin{theorem} \label{2ndmain}
Let $\underline{h} = (1,c,h_2,\dots,h_t,\dots, h_s)$ be an SI-sequence.  Then
the scheme $G$ described in Theorem \ref{1stmain} has maximal graded Betti
numbers among arithmetically Gorenstein subschemes of $\proj{n}$ with the
weak Lefschetz property and $h$-vector $\underline{h}$.
\end{theorem} 

If the embedding dimension satisfies $n+1 \geq s+c$ , we can choose the
Gorenstein scheme $G$ described in Theorem \ref{1stmain} as one defined by
a reduced monomial ideal. Thus $G$ corresponds to a simplicial complex. If
$P$ is a simplicial polytope then the Stanley-Reisner  ring  $K[\Delta(P)]$
of the boundary complex $\Delta(P)$ of $P$ is a reduced Gorenstein
$K$-algebra. The $h$-vector and the graded Betti numbers of $P$ are the
corresponding numbers of $K[\Delta(P)]$. The famous $g$-theorem mentioned
above characterizes the $h$-vectors of simplicial polytopes as
SI-sequences. For the Betti numbers we have the following result analogous
to Theorem \ref{2ndmain}. 

\begin{theorem} \label{3rdmain} Let $K$ be a field of characteristic zero
  and let $\underline{h}$ be an SI-sequence. Then there is a simplicial
  polytope with $h$-vector $\underline{h}$ having  maximal graded Betti
numbers among all simplicial polytopes with $h$-vector $\underline{h}$. 
\end{theorem} 

For the proof we have to relate our previous results to the work of Billera
and Lee in \cite{BL2}. It uses the combinatorial description of the
irreducible components of the arithmetically Gorenstein schemes occurring in
Theorem \ref{1stmain} (cf.\ Theorem \ref{maingorthm}). This description
allows an interpretation of a 
construction for simplicial polytopes with the help of Gorenstein linked
ideals. In the course of the proof of Theorem \ref{3rdmain} we also show
that the upper bounds on 
the Betti numbers for arbitrary graded  Cohen-Macaulay $K$-algebras with
given Hilbert function in
\cite{bigatti} and \cite{hulett} are best possible even for shellable
simplicial complexes. 

The maximal graded Betti numbers mentioned in Theorems \ref{2ndmain} and
\ref{3rdmain} can be computed effectively (cf.\ Theorems
\ref{main-thm-paper} and \ref{fax-9.5}). 

In spite of our results above and the conjectured generalization of the
$g$-theorem to simplicial spheres, we would like to pose the following
question: 

\begin{question}\label{reduced-WLP}
Does {\em every} reduced, arithmetically Gorenstein subscheme of projective
space possess the weak Lefschetz property?
\end{question} 
\medskip

\noindent
In case this question has an affirmative answer, this would imply that we have
characterized all Hilbert functions of reduced \aG schemes and that among
the  reduced
\aG schemes with fixed Hilbert function we always have the existence of a
scheme with maximal Betti 
numbers.
\smallskip

\begin{center} 
{\bf Acknowledgement} 
\end{center} 

The second author would like to thank J.\ Herzog for motivating 
discussions on simplicial polytopes. 


\section{Background}

We begin by recording the following notation and conventions.

\begin{itemize}

\item Let $K$ be any field.  We set

\begin{itemize}

\item[$\circ$] $R = K[x_0,\dots,x_n]$ and $\proj{n} = \Proj R$;

\item[$\circ$] $T = K[z_1,\dots,z_c]$ and  $\bar T = K[z_2,\dots,z_c]$, where
$c <n$ (this will be used primarily in section \ref{arb-hf-sect}).

\item[$\circ$] $S = K[y_1,\dots,y_c]$ (this will be used primarily in section
\ref{arb-hf-sect}).

\end{itemize}

\item Following \cite{KMMNP} we say that two homogeneous ideals $I_1$ and
$I_2$ are CI-linked (resp.\ G-linked) if they are linked using a complete
intersection (resp.\ a Gorenstein ideal which is not necessarily a complete
intersection).  That is, we require
\[
J \subset I_1 \cap I_2, \ \ \ \ J:I_1 = I_2, \ \ \ \ J:I_2 = I_1
\]
where $J$ is a complete intersection (resp.\ a Gorenstein ideal which is not
necessarily a complete intersection).

\item Let $V_1$ and $V_2$ respectively be the linked subschemes of $\proj{n}$
defined by the ideals $I_1$ and $I_2$ above.  Then $V_1$ and $V_2$ are also
said to be CI-linked (resp.\ G-linked).  If $X$ is the Gorenstein ideal
defined by $J$, we write $V_1 \stackrel{X}{\sim} V_2$.  The property of being
linked forces $V_1$ and $V_2$ to be equidimensional of the same dimension (as
$X$).  The link is {\em geometric} if $V_1$ and $V_2$ have no common
component.

\item If $X$ is a subscheme of $\proj{n}$ with saturated ideal $I_X$, and if
$t \in {\mathbb Z}$ then the Hilbert function of $X$ is denoted by
\[
h_X (t) = h_{R/I_X} (t) = \dim_K [R/I_X]_t .
\]

\item If $X$ is \acm of dimension $d$ then $A=R/I_X$ has Krull dimension
$d+1$ and a general set of $d+1$ linear forms forms a regular sequence for
$A$.   Taking the quotient of $A$ by such a regular sequence gives a
zero-dimensional Cohen-Macaulay ring called the {\em Artinian reduction} of
$A$ (or of $X$ -- cf.\ \cite{migliore}).  The Hilbert function of the Artinian
reduction of $R/I_X$ is called the {\em $h$-vector} of $R/I_X$ (or of $X$). 
This is a finite sequence of integers.

\item If $X$ is arithmetically Gorenstein with $h$-vector $(1,c,\dots,h_s)$
then this $h$-vector is symmetric ($h_s = 1$, $h_{s-1} = c$, etc.) and $s$ is
called the {\em socle degree} of $X$.

\end{itemize}

 The underlying idea of \cite{GM5} was to produce the desired
arithmetically Gorenstein scheme as a sum of CI-linked, arithmetically
Cohen-Macaulay codimension two schemes, chosen so that the complete
intersection linking them is already a union of linear subvarieties with nice
intersection properties.  This guarantees that the intersection of the linked
schemes will also be a good linear configuration.  The deformation result
arose as a consequence of work of Diesel \cite{diesel} about the
irreducibility of the Hilbert scheme of Gorenstein Artinian
$k$-algebras of codimension three with given Hilbert function and degrees of
generating sets.

Diesel's methods do not hold in higher codimension.  However, as a first step
in the above direction in higher codimension, we conjecture that every
Hilbert function which arises as the Hilbert function of an arithmetically
Gorenstein scheme in fact occurs for some arithmetically Gorenstein scheme
which is a generalized stick figure.  It is not yet known which Hilbert
functions in fact arise in this way, but a large class is known thanks to
\cite{harima} (as recalled below), and we will prove our conjecture for this
class.

Our approach is to combine ideas of \cite{harima} and \cite{GM5}, with
several new twists.  We recall some facts and definitions from the former.

\begin{definition}\label{basic-Gor-defs}
 Let $\underline{h} = (h_0,\dots,h_s)$ be a sequence of positive integers.
$\underline{h}$ is called a {\em Gorenstein sequence} if it is the Hilbert
function of some Gorenstein Artinian $k$-algebra.  $\underline{h}$ is {\em
unimodal} if $h_0 \leq h_1 \leq \dots \leq h_j \geq h_{j+1} \geq \dots \geq
h_s$ for some $j$.  $\underline{h}$ is called an SI-sequence (for
Stanley-Iarrobino) if it satisfies the following two conditions:

\begin{itemize}
\item[(i)] $\underline{h}$ is symmetric, i.e.\ $h_{s-i} = h_i$ for all
$i=0,\dots,\lfloor \frac{s}{2} \rfloor$.
\item[(ii)] $(h_0, h_1-h_0,h_2-h_1,\dots,h_j-h_{j-1})$ is an O-sequence,
where $j = \lfloor \frac{s}{2} \rfloor$; i.e.\ the ``first half'' of
$\underline{h}$ is a {\em differentiable} O-sequence.
\end{itemize}
\end{definition}

\begin{remark}\label{defoft}
 From the above definition of an SI-sequence it follows that
\[
(h_0, h_1-h_0,h_2-h_1,\dots,h_t-h_{t-1})
\]
 is an O-sequence, where
$t = \min \{ i | h_i \geq h_{i+1} \}$.  However, we cannot replace the given
condition by this one.  For instance, the sequence
\[
(1,3,6,6,7,6,6,3,1)
\]
satisfies this condition (with $t=2$) and is unimodal and symmetric, but the
last condition of the definition is not satisfied, and we do not want to allow
this kind of behavior.
\end{remark}

\begin{definition} \label{WSP}
A graded Artinian algebra $A = \bigoplus_{i-1}^s A_i$, with $A_s \neq 0$, is
said to have the {\em weak Lefschetz property} (WLP for short) if $A$ satisfies
the following two conditions.
\begin{itemize}
\item[(i)] The Hilbert function of $A$ is unimodal.
\item[(ii)] There exists $g \in A_1$ such that the $k$-vector space
homomorphism $g:A_i \rightarrow A_{i+1}$ defined by $f \mapsto gf$ is either
injective or surjective for all $i=0,1,\dots,s-1$.
\end{itemize}
Since we are interested in the question of lifting, we will say that {\em an
\acm subscheme $X$ of $\proj{n}$ has the  weak Lefschetz property} if there is
an Artinian reduction of $X$ having the weak Lefschetz property.
\end{definition}

\noindent In \cite{harima}, this property is also called the {\em
weak Stanley property}.  Harima proved the following very interesting result
(\cite{harima} Theorem 1.2).

\begin{theorem}  \label{harima-thm}
Let $\underline{h} = (h_0,h_1,\dots,h_s)$ be a sequence of positive integers.
Then $\underline{h}$ is the Hilbert function of some Gorenstein Artinian
$k$-algebra with the WLP if and only if $\underline{h}$ is an SI-sequence.
\end{theorem}

This provides a huge class of Hilbert functions for which it is known that
there is an Artinian Gorenstein ideal.  In codimension 3 it is known that all
Artinian Gorenstein ideals have Hilbert functions which are SI-sequences, and
in codimension $\geq 5$ it is known that not all do.  (Codimension 4 is open.)
Our goal is to provide, for each such SI-sequence in any codimension, a reduced
set of points (or more generally a reduced union of linear varieties) which is
arithmetically Gorenstein and whose $h$-vector is the given SI-sequence.  See
\cite{migliore} for more on $h$-vectors.

It is well known that the sum of the ideals of two geometrically linked,
arithmetically Cohen-Macaulay subschemes of $\proj{n}$ is arithmetically
Gorenstein of height one greater, whether they are CI-linked \cite{PS} or
G-linked (cf.\ \cite{migliore}).  Harima (\cite{harima}, Lemma 3.1) has
computed the Hilbert function of the Gorenstein ideals so obtained in the
case of CI-linkage; the proof is the same for G-linkage.  However, we would
like to record this result in a different (but equivalent) way, more in line
with our needs.

\begin{lemma}\label{comp-tool}
Let $V_1 \stackrel{X}{\sim} V_2$, where $X$ is arithmetically Gorenstein,
$V_1$ and $V_2$ are arithmetically Cohen-Macaulay of codimension $c$ with
saturated ideals $I_{V_1}$ and $I_{V_2}$, and the link is geometric.  Then
$I_{V_1} + I_{V_2}$ is the saturated ideal of an arithmetically Gorenstein
scheme $Y$ of codimension $c+1$.  The Hilbert functions are related as
follows.

Let $\underline{c} = (1,c,c_2,\dots,c_{s-1},c_s)$ be the $h$-vector of $X$;
note that $c_{s-1} = c$ and $c_s = 1$.  Let $\underline{g} =
(1,g_1,\dots,g_r)$ be the $h$-vector of $V_1$ (so we have $g_1 \leq c$,
with equality if and only if $V_1$ is non-degenerate, and $r \leq s$).  Let
$\underline{g}' = (1,g'_1,\dots)$ be the $h$-vector of $V_2$.  By \cite{DGO}
Theorem 3 (see also \cite{migliore} Corollary 5.2.19) we have the $h$-vector
of $V_2$ given by
\[
g'_i = c_{s-i} - g_{s-i}
\]
for $i \geq 0$.  Then the sequence $d_i = (g_i + g'_i - c_i)$ is the {\em
first difference} of the $h$-vector of $Y$.
\end{lemma}

\begin{example}\label{illustrate-comp}
A twisted cubic curve $V_1$ in $\proj{3}$ is linked to a line $V_2$ by the
complete intersection of two quadrics.  The intersection of these curves is
the arithmetically Gorenstein zeroscheme $Y$ consisting of two points.  This
is reflected in the diagram
\[
\begin{array}{ccccccccc}
X:    & 1 & 2 & 1 \\
V_1 : & 1 & 2   \\
V_2 : & 1 \\
\Delta Y :   & 1 & 0 & -1
\end{array}
\]
adding the second and third rows and subtracting the first to obtain the
fourth, and so the $h$-vector of $Y$ is (1,1), obtained by ``integrating''
the vector $(1,0 -1)$.  The notation $\Delta Y$ serves as a reminder that the
row is really the first difference of the $h$-vector of $Y$, and it will be
used in the remaining sections.
\end{example}

\begin{remark} \label{compare-harima-1}
Many papers on Hilbert functions or minimal free resolutions of Gorenstein
ideals (for instance \cite{GHS1}, \cite{GHS2} \cite{GKS}, \cite{GM5},
\cite{GPS}, \cite{GS}, \cite{harima2}) use as the method of constructing
Gorenstein ideals this notion of adding the ideals of geometrically CI-linked
Cohen-Macaulay ideals.  In codimension three this was enough, as shown in
\cite{GM5}.  However, in higher codimension this is not enough.  For
instance, in \cite{GKS} Remark 3.5 the authors say that with their
construction method they cannot obtain the $h$-vector
\[
1 \ \ 4 \ \ 10 \ \ n \ \ 10 \ \ 4 \ \ 1
\]
where $14 \leq n \leq 19.$  Indeed, one can check that taking $n=19$, it is
impossible to obtain such a Hilbert function as the sum of geometrically
CI-linked Cohen-Macaulay ideals (see Example \ref{ci-not-enough}).

Of course such an $h$-vector is a special case of
Harima's theorem.  Harima begins with the same complete intersection trick,
but he adds a very nice extra ingredient which we summarize briefly in Lemma
\ref{harima-lemma3.3}.  We would like to remark that to extend this approach
to the non-Artinian case, one would have to prove that the sum of CI-linked
ideals can be done in such a way that a very precisely determined subset of
the resulting Gorenstein set of points can always be found lying on a
hyperplane, and that the residual is again Gorenstein,  arriving at the
desired $h$-vector.  (Note, however, that Harima's trick works for a general
linear form, while in higher dimension the hyperplane is very special since
it contains a large number of the points.)  In this paper we do not take
this approach.  However, we show in Section \ref{subspace-sect} that the
configurations that we obtain have what we call the ``subspace property,''
and this corresponds precisely to the higher dimensional analog of Harima's
trick. See Remark \ref{compare-harima-2}.  This is a by-product of a
completely different approach using Gorenstein liaison.
\end{remark}

\begin{definition} Let $X$ be a finite set of points in $\proj{n}$ with
Hilbert function $h_X (i)$.  Then $\sigma(X) = \min \{ i | \Delta h_X (i) =
0 \}$.  
\end{definition}

\begin{lemma}[\cite{harima}, Lemma 3.3] \label{harima-lemma3.3}
Let $X$ and $Y$ be two finite sets of points in $\proj{n}$ such that $X \cap
Y = \emptyset$ and $X \cup Y$ is a complete intersection, and put $A =
K[x_0,\dots,x_n]/(I(X)+I(Y))$.  Furthermore put $a = \sigma(X)-1$, $b =
\sigma(X\cup Y) - \sigma(X)-1$, $c = \sigma(X \cup Y) - 1$.  Assume that $2
\sigma(X) \leq \sigma(X \cup Y)$ and $|X| \geq 2$.  Let $L \subset \proj{n}$
be a hyperplane defined by a polynomial $G \in R_1$ and let $g \in A_1$ be
the image of $G$.  Assume that $X \cap L = \emptyset$.  Let $d$ be an integer
such that $1 \leq d \leq \sigma(X \cup Y) - 2 \sigma(X)$ and let $[0:g^d]$
denote the homogeneous ideal generated by homogeneous elements $f \in A$ such
that $g^d f = 0$.  Then the Hilbert function $h_{A/[0:g^d]}(i)$ is a
Gorenstein SI-sequence as follows:
\begin{equation}\label{harima-hf2}
h_{A/[0:g^d]}(i) =
\left \{
\begin{array}{ll}
h_X(i), & \hbox{if $i = 0,1,\dots,a-1$,} \\
|X|, & \hbox{if $i = a,\dots,b-d$,} \\
h_X(c-1-i-d), & \hbox{if $i = b+1-d,\dots,c-1-d$,},
\end{array}
\right.
\end{equation}
and $\sigma(A/[0:g^d]) = c-d$.
\end{lemma}

\begin{remark}\label{def-reg} Let $X \subset \PP^n$ be a projective
subscheme. Then its  Castelnuovo-Mumford regularity is 
\[
\reg(X) := \reg(I_X) = \min\{ j \ | \ h^i (\proj{n}, {\mathcal I}_X (j-i)
= 0 \; \mbox{for all } i \; \mbox{with} \; 1 \leq i \leq \dim X +1 \ \}.
\]
If $X$ is arithmetically Cohen-Macaulay and its  $h$-vector  is 
$(1,h_1,\dots,h_s)$ with $h_s > 0$ then $\reg(X) = s+1$. Thus, we observe that 
$\sigma (X) = \reg(X)$ if $X$ is a zeroscheme.
\end{remark}


\section{Generalized stick figures and a useful construction}
\label{maximal-hf-sect}

In this section we give a construction which is useful for producing \acm
subschemes of projective space, especially unions of linear varieties with
``nice'' singularities.  It is an application of a result in \cite{KMMNP}.  As
a consequence, in the next section we show how to construct such unions of
linear varieties with ``maximal'' Hilbert function among the \acm schemes with
fixed regularity and initial degree.

The notion of a {\em stick figure} curve was introduced classically, and it
has culminated with the recent solution by Hartshorne of the so-called {\em
Zeuthen problem} \cite{hartshorne}.  A stick figure is simply a union of
lines, no three of which meet in a point.  We will make use of the following
type of configuration, which was introduced for codimension two in \cite{BM5}
and in this generality in \cite{MN2}.

\begin{definition} \label{gen-stick-fig}
Let $V$ be a union of linear subvarieties of $\proj{m}$ of the same dimension
$d$.  Then $V$ is a {\em generalized stick figure} if the intersection of any
three components of $V$ has dimension at most $d-2$ (where the empty set is
taken to have dimension $-1$).  In particular, if $d=1$ then $V$ is a stick
figure.
\end{definition}

\begin{remark}\label{can-add}
The property of being a generalized stick figure has a useful consequence in
liaison.  As mentioned in Lemma \ref{comp-tool}, if $V_1$ and $V_2$ are
arithmetically  Cohen-Macau\-lay subschemes of projective space which are
geometrically linked by an arithmetically Gorenstein scheme $X_c$ ($c$ is the
codimension) then $I_{V_1} + I_{V_2}$ is the saturated ideal of an
arithmetically Gorenstein subscheme $X_{c+1}$ of codimension $c+1$.  In this
paper we want to use this fact to construct many reduced arithmetically
Gorenstein subschemes of projective space.  A key observation for us is thus
the following.  If our ``linking scheme'' $X_c \subset \proj{n}$ is a
generalized stick figure, with $c<n$, then the sum of the linked ideals
defines the arithmetically Gorenstein scheme  $X_{c+1}$ which is again a
reduced union of linear varieties.  Indeed, the components of $X_{c+1}$ are
generically the intersection of two linear varieties.
\end{remark}

\begin{remark} \label{right-codim}
Let $\{L_1,\dots,L_p \}$ be a set of generically chosen linear forms in the
ring
$K[x_0,\dots,x_N]$ (i.e.\ hyperplanes in $\proj{N}$).  Let
\[
\begin{array}{rcl}
A_1 & = & (L_{i_1},\dots,L_{i_c}) \\
A_2 & = & (L_{j_1},\dots,L_{j_c}) \\
A_3 & = & (L_{k_1},\dots,L_{k_c})
\end{array}
\]
be three {\em different} ideals generated by subsets of $\{L_1,\dots,L_p \}$
defining codimension
$c$ linear varieties.  Then these three varieties meet in codimension $
c+1$ if and  only if there are exactly $c+1$ different linear forms among the
$3c$ generators. That is, the union of these three varieties violate the
condition of being a generalized stick figure if and only if there are exactly
$c+1$  different linear forms among the $3c$ generators of the ideals $A_1$,
$A_2$ and $A_3$.
\end{remark}

We now give our construction for \acm schemes.  We begin by recalling the
notion of {\em Basic Double G-linkage} introduced in \cite{KMMNP}, so called
because of part (iv) and the notion of Basic Double Linkage
(\cite{lazarsfeld-rao}, \cite{BM4}, \cite{GM4}).

\begin{lemma}[\cite{KMMNP} Lemma 4.8, Remark 4.9 and Proposition
5.10]\label{KMMNP-lemma} Let $J \supset I$ be homogeneous ideals of $R =
K[x_0,\dots,x_n]$, defining schemes $W \subset V \subset \proj{n}$ such that
$\codim V + 1 = \codim W$.  Let $A \in R$ be an element of degree $d$ such
that $I:A = I$.  Then we have
\begin{itemize}
\item[(i)] $\deg (I+A \cdot J) = d \cdot \deg I + \deg J$.

\item[(ii)] If $I$ is perfect and $J$ is unmixed then $I+A\cdot J$ is unmixed.

\item[(iii)] $J/I \cong [(I+A\cdot J)/I] (d).$

\item[(iv)] If $V$ is \acm with property $G_0$ and $J$ is unmixed then $J$
and $I + A\cdot J$ are linked in two steps using Gorenstein ideals.

\item[(v)] The Hilbert functions are related by
\[
\begin{array}{rcl}
h_{R/(I+A\cdot J)} (t) & = & h_{R/(I+(A))} (t) + h_{R/J}(t-d) \\
& = &h_{R/I}(t) - h_{R/I}(t-d) + h_{R/J}(t-d)
\end{array}
\]

\end{itemize}
\end{lemma}

Lemma \ref{KMMNP-lemma} should be interpreted as viewing the scheme $W$
defined by $J$ as a divisor on the scheme $V$ defined by $I$, and adding to it
a hypersurface section $H_A$ of $V$ defined by the polynomial $A$.  Note
that $I_{H_A} = I_V + (A)$.  If $V$ and $W$ are \acm then the divisor $W+H_A$
is again \acm (by step 4).  As an immediate application we have the following
by successively applying Lemma \ref{KMMNP-lemma}.

\begin{corollary} \label{hyperplane-sects}
Let $V_1 \subset V_2 \subset \dots \subset V_r \subset \proj{n}$ be \acm
schemes of the same dimension, each generically Gorenstein.  Let
$H_1,\dots,H_r$ be hypersurfaces, defined by forms $F_1,\dots,F_r$, such that
for each $i$,
$H_i$ contains no component  of $V_j$ for any $j \leq i$.  Let $W_i$ be the
\acm schemes defined by the corresponding hypersurface sections: $I_{W_i} =
I_{V_i} + (F_i)$.  Then

\begin{itemize}
\item[(i)] viewed as divisors on $V_r$, the sum $Z$ of the $W_i$ (which
is just the union if the hypersurfaces are general enough) is in the same
Gorenstein liaison class as $W_1$.

\item[(ii)] In particular, $Z$ is arithmetically Cohen-Macaulay.

\item[(iii)] As ideals we
have
\[
I_Z = I_{V_r} + F_r \cdot I_{V_{r-1}} + F_r F_{r-1} I_{V_{r-2}} + \cdots
+ F_r F_{r-1} \cdots F_2 I_{V_1} + (F_r F_{r-1}\cdots F_1).
\]

\item[(iv)] Let $d_i = \deg F_i$.  The Hilbert functions are related by the
formula
\[
\begin{array}{rcl}
h_Z (t) & = & h_{W_r} (t) + h_{W_{r-1}} (t-d_{r}) + h_{W_{r-2}}
(t-d_{r}-d_{r-1}) +
\dots \\
&& + h_{W_1} (t-d_{r}-d_{r-1}-\dots-d_2).
\end{array}
\]
\end{itemize}
\end{corollary}

\begin{remark}
Parts (ii), (iii) and (iv) of Corollary \ref{hyperplane-sects} have been proved
independently by Ragusa and Zappal\`a (\cite{RZ} Lemma 1.5).  
\end{remark}

\newsavebox{\build}
\savebox{\build}(400,85)[tl]
{
\begin{picture}(400,125)
\put (120,20){\line (0,1){60}}
\put (140,20){\line (0,1){60}}
\put (160,20){\line (0,1){60}}
\put (180,20){\line (0,1){60}}
\put (200,20){\line (0,1){60}}
\put (220,20){\line (0,1){60}}
\put (240,20){\line (0,1){60}}
\put (260,20){\line (0,1){60}}
\put (280,20){\line (0,1){60}}

\put (110,30){\line (1,0){180}}
\put (110,40){\line (1,0){180}}
\put (110,50){\line (1,0){180}}
\put (110,60){\line (1,0){180}}
\put (110,70){\line (1,0){180}}

\put (117,27){$\bullet$}
\put (117,37){$\bullet$}
\put (117,47){$\bullet$}
\put (117,57){$\bullet$}
\put (117,67){$\bullet$}
\put (137,27){$\bullet$}
\put (137,37){$\bullet$}
\put (137,47){$\bullet$}
\put (137,57){$\bullet$}
\put (157,27){$\bullet$}
\put (157,37){$\bullet$}
\put (157,47){$\bullet$}
\put (157,57){$\bullet$}
\put (177,27){$\bullet$}
\put (177,37){$\bullet$}
\put (197,27){$\bullet$}
\put (197,37){$\bullet$}
\put (217,27){$\bullet$}
\put (217,37){$\bullet$}
\put (237,27){$\bullet$}
\put (257,27){$\bullet$}

\put (92,27){$\scriptstyle L_1$}
\put (92,37){$\scriptstyle L_2$}
\put (97,50){$\scriptscriptstyle \vdots$}
\put (92,67){$\scriptstyle L_d$}

\put (112,7){$\scriptstyle M_1$}
\put (132,7){$\scriptstyle M_2$}
\put (152,7){$\scriptstyle M_3$}
\put (172,7){$\scriptstyle M_4$}
\put (192,7){$\scriptstyle M_5$}
\put (212,7){$\scriptstyle M_6$}
\put (232,7){$\scriptstyle M_7$}
\put (257,7){$\scriptstyle \dots$}
\put (272,7){$\scriptstyle M_e$}
\end{picture}
}

\begin{corollary}\label{construct-acm}
Let $R = K[x_0,\dots,x_n]$.  Consider the complete intersection $(A,B)$ in
$R$, where $A = \prod_{i=1}^d L_i$ and $B = \prod_{i=1}^e M_i$.
Thinking of the $L_i$ and $M_i$ as hyperplanes which are pairwise linearly
independent, the intersection of any
$L_i$ with $M_k$ is a codimension two linear variety, $P_{i,k}$.
Consider a union $Z$ of such varieties subject to the condition that if
$P_{i,k} \subset Z$ then $P_{j,\ell} \subset Z$ for all $(j,\ell)$ satisfying
$j \leq i$ and $\ell \leq k$:

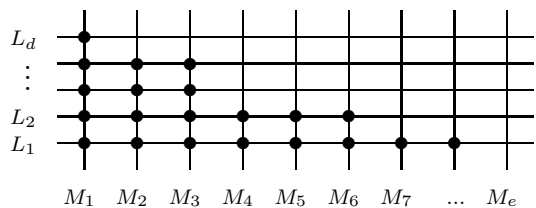
\begin{picture}(400,115)
\put(20,50){\usebox{\build}}
\end{picture}

\noindent Then $Z$ is arithmetically Cohen-Macaulay.
\end{corollary}

\begin{proof}
Apply Corollary \ref{hyperplane-sects}, using $V_1 = M_1$, $V_2 = M_1 \cup
M_2$, etc.\ and taking the hypersurface sections obtained by suitable
multiples of the $L_i$ (starting with $L_d$ and working backwards), working
our way down the picture.   This corollary can also be obtained using the
lifting results of \cite{MN2}.
\end{proof}

\begin{remark}
In the case $n=2$ these are not necessarily the $k$-configurations of
\cite{harima2}, since consecutive lines $L_i$ are allowed to have the same
number of points.  This situation has been studied in \cite{GPS}, and such a
configuration in $\proj{2}$ is a special case of a {\em weak
$k$-configuration}.  (In general a weak configuration does not require the
existence of the ``vertical'' $M_i$.)  It was shown in \cite{GPS} that the
Hilbert function of a weak $k$-configuration in $\proj{2}$ can be immediately
read from the number of points on each line.  This type of special case was
also extended to higher codimension in \cite{MN2}, section 4.

In the case of lines in $\proj{3}$, Corollary \ref{construct-acm} is
essentially contained in \cite{giorgio}.
\end{remark}


\section{A construction of \acm and Gorenstein ideals with ``maximal''
Hilbert function} \label{max-hf-sect}

In this section we show how to construct certain Gorenstein subschemes of
projective space with two properties: they are generalized stick figures and
they have very special Hilbert function (essentially maximal until halfway
through the $h$-vector).  The Hilbert function that we seek, for codimension
$c$, corresponds to the $h$-vector
\begin{equation}\label{max-Gor-hf}
\left ( 1, c, \binom{c+1}{c-1}, \dots,
 \binom{c-1+t}{c-1},\dots, \binom{c-1+t}{c-1} , \dots, \binom{c+1}{c-1},c,1
\right )
\end{equation}
where the terms in the middle are all equal to $\binom{c-1+t}{c-1}$ and the
last non-zero entry is in degree $s$.  Then we have
that the first occurrence of the value  $\binom{c-1+t}{c-1}$ is in degree $t$
and the last in degree $s-t$.

In order to construct such an ideal, we will first give a construction of
\acm ideals with maximal Hilbert function (expressed as an $h$-vector)
\begin{equation}\label{max-h-vtr}
\left ( 1, c, \binom{c+1}{2}, \dots, \binom{c-1+t}{t} \right ).
\end{equation}
The procedure is inductive, producing first a suitable answer in codimension
one, then codimension two, then codimension three, etc.  We will construct, in
codimension $c$, an \acm generalized stick figure $Z_{c,t}$ with ``maximal''
$h$-vector given by (\ref{max-h-vtr}).

\begin{theorem} \label{acm-scheme}
Let $R$ be a polynomial ring of dimension $n+1>c$.  Choose a set of $2t+c$
linear forms
\[
{\mathcal M}_{c,t} = \left \{ M_0,\dots,M_{t+\lfloor \frac{c-1}{2} \rfloor}
,L_0,\dots,L_{t+\lfloor
\frac{c-2}{2} \rfloor} \right \}
\]
in $R$.  Define
\begin{equation} \label{L<M}
\begin{array}{rll}
I_{Z_{c,t}} & = \displaystyle \bigcap_{0 \leq i_1
\leq i_2 < i_3 \leq i_4 < \dots < i_{c-1} \leq i_c \leq t+ \frac{c-2}{2}}
(M_{i_1}, L_{i_2}, M_{i_3}, L_{i_4}, \dots, L_{i_c}) & \hbox{
if $c$ is even}  ,\\
&& \\
I_{Z_{c,t}} & = \displaystyle \bigcap_{0 \leq i_1
\leq i_2 < i_3 \leq i_4 < \dots \leq i_{c-1} < i_c \leq t+ \frac{c-1}{2} }
(M_{i_1}, L_{i_2}, M_{i_3}, L_{i_4}, \dots, M_{i_c}) &
\hbox{ if  $c$ is odd.}
\end{array}
\end{equation}
(If $c=1$ the above range is understood to be $0 \leq i_1 \leq t$.)  If each
subset consisting of
$c+1$ of the
$2t+c$ linear forms generates an ideal of codimension $c+1$ then
$I_{Z_{c,t}}$ is a reduced Cohen-Macaulay ideal with $h$-vector
\[
\left ( 1, c, \binom{c+1}{2}, \dots, \binom{c+t-1}{t} \right ).
\]
If each subset consisting of $c+2$ of the $2t+c$ linear
forms generates a complete intersection of codimension $c+2$ then
$I_{Z_{c,t}}$ defines a generalized stick figure.

\end{theorem}

\begin{proof}

Consider the complete intersection $(A,B)$ in codimension 2, where
\[
A = \prod_{i=0}^{t+\lfloor \frac{c-1}{2} \rfloor} M_i \hbox{\ \ and \ \ } B =
\prod_{i=0}^{t+\lfloor \frac{c-2}{2} \rfloor} L_i
\]
(note that both products start with $i=0$).  Let us denote the scheme
defined by $(A,B)$ by $G_2$ (subscripts here will refer to the codimension).

When $c=2$, the scheme $Z_{2,t}$ has the following form (where the components
of $Z_{2,t}$ are represented by dots, and the set of all intersection points
is $G_2$):

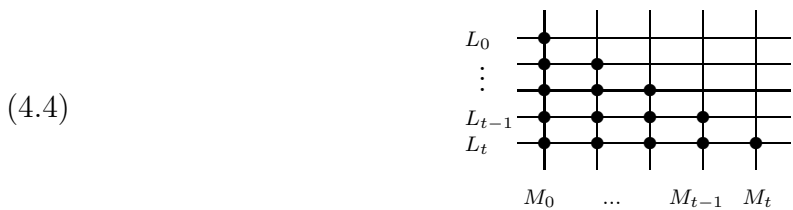
\newsavebox{\buildnice}
\savebox{\buildnice}(400,-20)[tl]
{
\begin{picture}(400,10)
\put (120,20){\line (0,1){60}}
\put (140,20){\line (0,1){60}}
\put (160,20){\line (0,1){60}}
\put (180,20){\line (0,1){60}}
\put (200,20){\line (0,1){60}}

\put (110,30){\line (1,0){103}}
\put (110,40){\line (1,0){103}}
\put (110,50){\line (1,0){103}}
\put (110,60){\line (1,0){103}}
\put (110,70){\line (1,0){103}}

\put (117,27){$\bullet$}
\put (117,37){$\bullet$}
\put (117,47){$\bullet$}
\put (117,57){$\bullet$}
\put (117,67){$\bullet$}
\put (137,27){$\bullet$}
\put (137,37){$\bullet$}
\put (137,47){$\bullet$}
\put (137,57){$\bullet$}
\put (157,27){$\bullet$}
\put (157,37){$\bullet$}
\put (157,47){$\bullet$}
\put (177,27){$\bullet$}
\put (177,37){$\bullet$}
\put (197,27){$\bullet$}

\put (90,27){$\scriptstyle L_t$}
\put (90,37){$\scriptstyle L_{t-1}$}
\put (95,50){$\scriptscriptstyle \vdots$}
\put (90,67){$\scriptstyle L_{0}$}

\put (112,7){$\scriptstyle M_0$}
\put (142,7){$\scriptstyle \dots$}
\put (167,7){$\scriptstyle M_{t-1}$}
\put (195,7){$\scriptstyle M_{t}$}
\end{picture}
}

\begin{equation} \label{pict-in-first-thm}
\begin{picture}(400,40)
\put(40,-10){\usebox{\buildnice}}
\end{picture}
\end{equation}

\vskip .5in

\noindent Using Corollary \ref{construct-acm}, it follows immediately that the
$h$-vector of $Z_{2,t}$ is
\[
\left (
1,2,3,4,\dots,t,t+1
\right )
\]
as claimed.

Now consider the case  $c=3$.  Let $V_1$ be the subset of $Z_{2,t}$ consisting
of all ``dots'' lying on $L_0$ (there is just one).  Let $V_2$ be the union of
the components of $X$  lying on {\em either} $L_0$ {\em or} $L_{1}$, and so
on.  Clearly we have
\[
V_1 \subset V_2 \subset \cdots \subset V_{t+1}
\]
and all $V_i$ are arithmetically Cohen-Macaulay.  We will apply Corollary
\ref{hyperplane-sects}, so to that end let us set $F_i = M_{i}$ ($1 \leq i
\leq t+1$) and let $W_i$ the hyperplane section of $V_i$ with $F_i$.  Then we
have
\[
\begin{array}{rcl}
\deg W_1 & = & 1 \\
\deg W_2 & = & 1+2=3 \\
& \vdots \\
\deg W_{t+1}& = & 1+2+\dots+(t+1) = \binom{t+2}{2}
\end{array}
\]
and by Corollary \ref{hyperplane-sects} the union of the $W_i$ is \acm of
codimension 3 with the $h$-vector
\[
\left (
1 , 3 , \dots, \binom{t+2}{2}
\right )
\]
 (use Lemma \ref{KMMNP-lemma}).
Clearly we also have
\[
W_1 \subset W_1 \cup W_2 \subset W_1 \cup W_2 \cup W_3 \cdots.
\]
Remark \ref{right-codim} quickly shows that each of these is a
generalized stick figure.  But clearly the union of the $W_i$ has the form
described in (\ref{L<M}), so this is $Z_{3,t}$ and we have completed the case
$c=3$.

To pass to codimension 4 using Corollary \ref{hyperplane-sects} again, we now
take $V_1,\dots,V_t$ to be the codimension 3 schemes just produced (i.e.\
set $V_1$ now to be the $W_1$ just produced, $V_2$ to be the $W_1 \cup W_2$
just produced, etc.), and we take $F_i = L_i$ ($1\leq i \leq t+1$).   As
before we have a generalized stick figure in codimension 4 which is \acm and
has maximal Hilbert function.   Note that the components of the \acm scheme
$Z_{4,t}$ we have produced have precisely the form described in (\ref{L<M}).

We continue by induction to finally produce the desired \acm scheme of
codimension $c$ with $h$-vector given by (\ref{max-h-vtr}).  (Note that the
range for $F_i$ changes: for instance, in codimension 5 we have $F_i =
M_{i}$ ($2 \leq i \leq t+2$).  It is a generalized stick figure by Remark
\ref{right-codim}.
\end{proof}

\begin{remark}   \label{codim-c-1-to-codim-c}
We will now turn to the construction of the arithmetically Gorenstein
generalized stick figure with $h$-vector given by (\ref{max-Gor-hf}).  First
we check numerically what is needed, using Lemma \ref{comp-tool} (see also
Example \ref{illustrate-comp}).  In the table below, the values occurring
between degrees $t+1$ and $s-t$ on any row are constant.  $G_{c-1}$ is the
$h$-vector of a codimension $c-1$ arithmetically Gorenstein scheme which
links \acm schemes $Z_{c-1}$ to $Y_{c-1}$, and the sum of linked ideals gives
an arithmetically Gorenstein scheme $G_c$ with the $h$-vector claimed in
(\ref{max-Gor-hf}).
\smallskip

\begin{tabular}{r|ccccccccccccc}
&\multicolumn{3}{l}{degree:}  \\
 &$\scriptstyle 0$ &$\scriptstyle 1$ &
$\scriptstyle 2$ & $\scriptstyle
\dots$ &
$\scriptstyle t$ &$\scriptstyle t+1$ &
$\scriptstyle \dots$ & $\scriptstyle s-t$ &$\scriptstyle s-t+1$ &$\scriptstyle
\dots$ &$\scriptstyle s$ &$\scriptstyle s+1$
\\ \hline
$G_{c-1}$  & $\scriptstyle 1$ & $\scriptstyle c-1$ &
$\binom{c}{c-2}$ &
$\scriptstyle \dots$ & $\binom{c-2+t}{c-2}$ & $\binom{c-2+t}{c-2}$ &
$\scriptstyle
\dots$ & $\binom{c-2+t}{c-2}$ & $\binom{c-2+t}{c-2}$ & $\scriptstyle \dots$ &
$\scriptstyle c-1$ & $\scriptstyle 1$ \\
$Z_{c-1}$ & $\scriptstyle 1$ & $\scriptstyle c-1$ & $\binom{c}{c-2}$ &
$\scriptstyle \dots$ & $\binom{c-2+t}{c-2}$ \\
$Y_{c-1}$ & $\scriptstyle 1$ & $\scriptstyle c-1$ & $\binom{c}{c-2}$ &
$\scriptstyle \dots$ & $\binom{c-2+t}{c-2}$ & $\binom{c-1+t}{c-2}$ &
$\scriptstyle
\dots$ & $\binom{c-1+t}{c-2}$ \\
$\Delta G_{c}$ &  $\scriptstyle 1$ & $\scriptstyle c-1$ & $\binom{c}{c-2}$
&
$\scriptstyle \dots$ & $\binom{c-2+t}{c-2}$ & $\scriptstyle 0$ & $\scriptstyle
\dots$ & $\scriptstyle 0$ & $- \binom{c-2+t}{c-2}$ & $\scriptstyle \dots$ &
$\scriptstyle -(c-1)$ & $\scriptstyle -1$ \\
$G_{c}$ &  $\scriptstyle 1$ & $\scriptstyle c$ & $\binom{c+1}{c-1}$
&
$\scriptstyle \dots$ & $\binom{c-1+t}{c-1}$ & $\binom{c-1+t}{c-1}$ &
$\scriptstyle
\dots$ & $\binom{c-1+t}{c-1}$ & $ \binom{c-2+t}{c-1}$ & $\scriptstyle \dots$
&
$\scriptstyle 1$ & $\scriptstyle 0$
\end{tabular}
\smallskip

\noindent 
Note that in the above table, if we take as
$G_{c-1}$ a Gorenstein scheme which ``levels off'' in degree $t+1$ rather than
degree $t$, then the residual Hilbert function changes, but the
Hilbert function of the sum of the linked ideals, $G_c$, does not change. 

The basic idea of our proof will be to assume by induction that a
generalized stick figure $G_{c-1}$ exists with the desired Hilbert
function, and that it contains $Z_{c-1,t}$, hence giving us a geometric link.
Then we have to show that the scheme $G_{c}$ obtained by adding the linked
ideals is again a generalized stick figure and contains $Z_{c,t}$ (which we
have already described).  This allows for the construction in the next
codimension, hence completing the induction.  In fact, we are able to give the
components of $G_{c}$ very explicitly!  We will refine the notation for $G_c$
in the proof, to account for more data.
\end{remark}

\begin{theorem} \label{max-gor}
Let $R$ be a polynomial ring of dimension $n+1>c \geq 1$.  Let $s \geq 2t$. 
Let
\[
{\mathcal N}_{c,s,t} = \left \{ M_0,\dots,M_{t + \lfloor \frac{c-1}{2}
\rfloor},L_0,\dots,L_{s-t+\lfloor \frac{c-2}{2} \rfloor} \right \}
\]
be a subset of $s+c$ linear forms in $R$.  Define an ideal $I_{G_{c,s,t}}$ as
follows.  If $c$ is even then $I_{G_{c,s,t}} = A_{c,s,t} \cap
B_{c,s,t} \cap C_{c,s,t}$ where
\[
\begin{array}{rcl}
A_{c,s,t}& = & \displaystyle \bigcap_{0 \leq i_1 \leq i_2 < i_3 \leq \cdots <
i_{c-1} \leq i_c \leq t+\frac{c-2}{2}} (M_{i_1},L_{i_2},\dots,M_{i_{c-1}},
L_{i_c}) , \\  \\
B_{c,s,t} & = & \displaystyle \bigcap_{0 \leq i_1 < i_2 \leq i_3 < \cdots \leq
i_{c-1}< i_c \leq t+\frac{c-2}{2}} (L_{i_1},M_{i_2},\dots,L_{i_{c-1}},
M_{i_c}), \\ \\
C_{c,s,t} & = & \displaystyle \bigcap_{
\begin{array}{c}
\scriptstyle 0 \leq i_1 \leq i_2 < i_3 \leq \cdots < i_{c-1} \leq
t+\frac{c-2}{2} \\
\scriptstyle t+\frac{c}{2} \leq i_c \leq s-t+\frac{c-2}{2}
\end{array}
}
(M_{i_1},L_{i_2},\dots,M_{i_{c-1}}, L_{i_c}).
\end{array}
\]
If $c$ is odd then $I_{G_{c,s,t}} = A_{c,s,t}' \cap
B_{c,s,t}' \cap C_{c,s,t}'$ where
\[
\begin{array}{rcl}
A_{c,s,t}'& = & \displaystyle \bigcap_{0 \leq i_1 \leq i_2 < i_3 \leq \cdots
\leq i_{c-1} < i_c \leq t+\frac{c-1}{2}}
(M_{i_1},L_{i_2},\dots,L_{i_{c-1}}, M_{i_c}) , \\  \\
B_{c,s,t}' & = & \displaystyle \bigcap_{0 \leq i_1 < i_2 \leq i_3 < \cdots
< i_{c-1}\leq i_c \leq t+\frac{c-3}{2}} (L_{i_1},M_{i_2},\dots,M_{i_{c-1}},
L_{i_c}), \\ \\
C_{c,s,t}' & = & \displaystyle \bigcap_{
\begin{array}{c}
\scriptstyle 0 \leq i_1 \leq i_2 < i_3 \leq \cdots \leq i_{c-1} \leq
t+\frac{c-3}{2} \\
\scriptstyle t+\frac{c-1}{2} \leq i_c \leq s-t+\frac{c-3}{2}
\end{array}
}
(M_{i_1},L_{i_2},\dots,L_{i_{c-1}}, L_{i_c}).
\end{array}
\]
If each subset of ${\mathcal N}_{c,s,t}$ consisting of $c+1$ elements
generates a complete intersection then $I_{G_{c,s,t}}$ is a reduced
Gorenstein ideal contained in
$I_{Z_{c,t}}$ with $h$-vector
\[
\left ( 1,c,\binom{c+1}{2}, \dots \right., \underbrace{\binom{c+t-1}{t},
\dots,\binom{c+t-1}{t}}_{\hbox{\rm flat part}}, \left. \dots,
\binom{c+1}{2},c,1
\right )
\]
where the final \hbox{\rm ``$1$''} occurs in degree $s$,  and $t+1$ is the
initial degree of the ideal if $c \geq 2$.  If  each subset of ${\mathcal
N}_{c,s,t}$ consisting of $c+2$ elements generates a complete intersection
of codimension
$c+2$ then $G_{c,s,t}$ is a generalized stick figure.
\end{theorem}

\begin{proof}
It is clear from the description that $I_{G_{c,s,t}}$ is a reduced ideal. 

Next we verify that if ${\mathcal N}_{c,s,t}$ has the property that each
subset of $c+2$ elements generates a complete intersection of codimension
$c+2$ then $G_{c,s,t}$ is a generalized stick figure.  To prove this, consider
first the case $c$ even.  Suppose we  have three components 
\[
\begin{array}{c}
P = (p_1,p_2,\dots,p_c) \\
Q = (q_1 ,q_2,\dots,q_c) \\
R = (r_1,r_2,\dots,r_c)
\end{array}
\]
(each taken from $A_{c,s,t}$, $B_{c,s,t}$ or $C_{c,s,t}$).  In order for
$P$ and $Q$ to have $c-1$ entries in common, we must be able to take out an
entry from $P$ and replace it with a  new entry, giving $Q$.  Putting both of
these entries in and removing a third one must give $R$.  Because of the
rigid form of the components given in the statement of the theorem, one can
just check that this is impossible.

Next suppose $c$ is odd and $s=2t$.  Note that in this case $C'_{c,s,t}$ is
empty.  Then the argument is similar to the one above.

The case where $c$ is odd and $s>2t$ is handled similarly.  Here is it
slightly tricky to prove that we get a generalized stick figure because of
the fact that condition $C'_{c,s,t}$ allows two consecutive $L's$ at the
end.  But the subscript of the second one is bounded below by $t +
\frac{c-1}{2}$, and this fact is needed to complete the proof.  For example,
the linear forms \label{obst-to-gen-sf}
\[
\begin{array}{c}
(M_1,L_3,M_5,L_7,M_8) \\
(M_1,L_3,M_5,L_7,L_9) \\
(L_3,M_5,L_7,M_8,L_9)
\end{array}
\]
seem at first glance to be a counter-example, but the ``9" in $L_9$ cannot
simultaneously be bounded above by $t+\frac{c-3}{2}$ and below by
$t+\frac{c-1}{2}$.

The fact that $Z_{c,t} \subset G_{c,s,t}$ follows from the observation that
$A_{c,s,t}$ (resp.\ $A'_{c,s,t}$) is $I_{Z_{c,t}}$, thanks to Theorem
\ref{acm-scheme}.

It remains to show that $I_{G_{c,s,t}}$ is a Gorenstein ideal with the
correct $h$-vector.  For this we proceed by induction on $c \geq 1$.  If
$c=1$ then $I_{G_{1,s,t}} = (M_0 \cdot \hdots \cdot M_t \cdot L_0 \cdot
\hdots \cdot L_{s-t-1})$ is a principal ideal of degree $s+1$, thus having
the claimed properties.

If $c=2$ then $I_{G_{2,s,t}} = (M_0 \cdot \hdots \cdot M_t,L_0 \cdot \hdots
\cdot L_{s-t})$ is a complete intersection by assumption on ${\mathcal
N}_{c,s,t}$.  It is easy to check its $h$-vector.

Now suppose that $c \geq 3$.  We distinguish two cases.

\underline{Case 1}:
Assume that $c$ is odd.  We have at our disposal the set
\[
{\mathcal N}_{c,s,t} = \left \{
M_0,\hdots,M_{t+ \frac{c-1}{2}} ,L_0,\hdots,L_{s-t+\frac{c-3}{2}}
\right \}
\]
where each subset of $c+1$ elements is linearly independent (over $K$).  Let
us temporarily re-name the linear form $M_{t+\frac{c-1}{2}}$ with the new
name $L_{s-t+\frac{c-1}{2}}$.  Hence we now have the set
\[
{\mathcal N}_{c-1,s+1,t} = \left \{
M_0,\hdots,M_{t+\frac{c-3}{2}},L_0,\hdots,L_{s-t+\frac{c-1}{2}}
\right \}.
\]
By induction we can use ${\mathcal N}_{c-1,s+1,t}$ to get the Gorenstein
ideal $I_{G_{c-1,s+1,t}}$ which defines a generalized stick figure.  The
configuration $Z_{c-1,t}$ is formed using the set ${\mathcal M}_{c-1,t} = \{
M_0,\hdots,M_{t+\frac{c-3}{2}} L_0,\hdots,L_{t+\frac{c-3}{2}} \}$.  Since $s
\geq 2t$, we observe that
\begin{equation}\label{stillincl}
{\mathcal M}_{c-1,t} \subset {\mathcal N}_{c-1,s+1,t} \backslash \{
L_{s-t+\frac{c-1}{2}} \}.
\end{equation}
Now we re-name $L_{s-t+\frac{c-1}{2}}$ back to $M_{t+\frac{c-1}{2}}$. 
The configuration $G_{c-1,s+1,t}$ then becomes (with the re-naming) 
a configuration, which we now call $G'_{c-1,s+1,t}$, whose ideal
is 
\begin{equation}\label{listcomps}
I_{G'_{c-1,s+1,t}} = A_{c-1,s+1,t} \cap B_{c-1,s+1,t} \cap \tilde
C_{c-1,s+1,t} \cap \tilde D_{c-1,s+1,t}
\end{equation}
where
\[
\begin{array}{rcl}
\tilde C_{c-1,s+1,t}&  =  &
\displaystyle 
\bigcap_{
\begin{array}{c}
\scriptstyle 0 \leq i_1 \leq i_2 < i_3 \leq \cdots < i_{c-2} \leq
t+\frac{c-3}{2} \\
\scriptstyle t+\frac{c-1}{2} \leq i_{c-1} \leq s-t+\frac{c-3}{2}
\end{array}
}
(M_{i_1},L_{i_2},\dots,M_{i_{c-2}}, L_{i_{c-1}}) \\ \\
\tilde D_{c-1,s+1,t}&  =  &
\displaystyle 
\bigcap_{
\begin{array}{c}
\scriptstyle 0 \leq i_1 \leq i_2 < i_3 \leq \cdots < i_{c-2} \leq
t+\frac{c-3}{2}\end{array}
}
(M_{i_1},L_{i_2},\dots,M_{i_{c-2}}, M_{t+\frac{c-1}{2}}),
\end{array}
\]
where it is understood that $\tilde C_{c-1,s+1,t} = R$ if $s=2t$.

 By the observation (\ref{stillincl}), this re-naming does not affect any
component of $Z_{c-1,t}$, which is defined by $A_{c-1,s+1,t}$.  Of course the
naming of the linear forms does not affect the properties of the
configurations, and so all together we have  that $G'_{c-1,s+1,t}$ is an
arithmetically Gorenstein generalized stick figure with the same $h$-vector
as $G_{c-1,s+1,t}$ and containing $Z_{c-1,t}$.

Since $G'_{c-1,s+1,t}$ is a generalized stick figure containing $Z_{c-1,t}$,
it provides a geometric link to the residual
\[
Y = G'_{c-1,s+1,t} \backslash Z_{c-1,t}.
\]
By Remark \ref{codim-c-1-to-codim-c}, the ideal $I_{G_{c,s,t}}$ given in the
statement of the theorem will have the desired properties if we can show that
$I_{G_{c,s,t}} = I_{Z_{c-1,t}} + I_Y.$

We have to show
\begin{itemize}
\item[(a)] For every choice of a component from $A_{c-1,s+1,t}$ (i.e.\ from
$Z_{c-1,t}$) and a component from either $B_{c-1,s+1,t}$, $\tilde
C_{c-1,s+1,t}$ or $\tilde D_{c-1,s+1,t}$ {\em which meet in codimension $c$},
their intersection occurs in either $A'_{c,s,t}$, $B'_{c,s,t}$ or
$C'_{c,s,t}$.

\item[(b)] Every component of $A'_{c,s,t}$, $B'_{c,s,t}$, or $C'_{c,s,t}$ is
the intersection of a component from $A_{c-1,s+1,t}$ and a component from
either $B_{c-1,s+1,t}$, $\tilde C_{c-1,s+1,t}$ or $\tilde D_{c-1,s+1,t}$.
\end{itemize}

All of these involve an analysis of how many linear forms can be common to
two (or three) of the components of our configurations.  For (a), in order for
the two components in codimension $c-1$ to meet in codimension $c$, their
intersection must be defined by only $c$ linear forms, so they must have
exactly $c-2$ linear forms in common.  We have to determine all the ways that
this can happen and show that we always get a component of $G_{c,s,t}$.  This
analysis then works backwards to show just how to produce any component of
$G_{c,s,t}$ as the intersection of two components of $G'_{c-1,s+1,t}$, which
answers (b).  

To answer these two questions, we leave it to the reader to verify that in
order to meet in codimension~$c$,

\begin{itemize}
\item a component from $A_{c-1,s+1,t}$ (i.e.\ $Z_{c-1,t}$) meets a component
 from $B_{c-1,s+1,t}$ either in the form $B'_{c,s,t}$ (and {\em all}
components of $B'_{c,s,t}$ arise in this way) or in the form $A'_{c,s,t}$
(and all components of $A'_{c,s,t}$ arise in this way {\em except} those
where $M_{i_c} = M_{t+\frac{c-1}{2}}$).

\item a component from $A_{c-1,s+1,t}$ meets a component from
$\tilde C_{c-1,s+1,t}$ in the form $C'_{c,s,t}$ (and {\em all} components of
$C'_{c,s,t}$ arise in this way).

\item a component from $A_{c-1,s+1,t}$ meets a component from $\tilde
D_{c-1,s+1,t}$ in the form $A'_{c,s,t}$ where $M_{i_c} =
M_{t+\frac{c-1}{2}}$, taking care of those components which were ``missing''
 from the first set above.
\end{itemize}

\noindent To do this one checks how it is possible to remove one entry from
the first component and one entry from the second and have
the remaining entries equal. (There are very few possibilities.) 

\underline{Case 2}:
Assume that $c$ is even.  Then ${\mathcal N}_{c,s,t} = {\mathcal
N}_{c-1,s+1,t}$.  By induction, the arithmetically Gorenstein scheme
$G_{c-1,s+1,t}$ contains $Z_{c-1,t}$.  We put $Y := G_{c-1,s+1,t} \backslash
Z_{c-1,t}$, the residual scheme.  The assertion follows because
$I_{G_{c,s,t}} = I_{Z_{c-1,t}} + I_Y$, which can be shown as in Case 1 (and
is easier).

This completes the construction of the arithmetically Gorenstein generalized
stick figure with ``maximal'' $h$-vector.
\end{proof}

\begin{example}
To construct the $h$-vector $(1,4,4,1)$ we take $s=3$,
$t=1$,  $c=4$,  and we get the components
$(M_0,L_0,M_1,L_1)$, $(M_0,L_0,M_1,L_2)$, $(M_0,L_0,M_2,L_2)$,
$(M_0,L_1,M_2,L_2)$, $(M_1,L_1,M_2,L_2)$,
$(L_0,M_1,L_1,M_2)$,
$(M_0,L_0,M_1,L_3)$,
$(M_0,L_0,M_2,L_3)$, $(M_0,L_1,M_2,L_3)$, $(M_1,L_1,M_2,L_3)$.
\end{example}

\begin{remark}
As was the case with Theorem \ref{acm-scheme}, the construction of the \aG
scheme with ``maximal Hilbert function'' described in Theorem \ref{max-gor}
can be viewed in a very concrete, geometrical way, especially for low
codimension.  For example, let us produce $G_{3,s,t}$.  We start with the set
\[
{\mathcal N}_{3,s,t} = \{ M_0,\dots,M_{t+1},L_0,\dots,L_{s-t} \}
\]
The result of renaming $M_{t+1}$ to $L_{s-t+1}$, considering $G_{2,s+1,t}$
and renaming back is that we have the complete intersection $(A,B)$ where 
\[
A = \prod_{i=0}^{t} M_i , \ \ \
B = \left ( \prod_{i=0}^{s-t} L_i \right ) \cdot M_{t+1}.
\]
The scheme $Z_{2,t}$ is a subconfiguration.  Let $Y$ be the
residual to $Z_{2,t}$ in this complete intersection.  In the following
diagram, $Z_{2,t}$ is given by the dots and $Y$ by the
intersection points without dots:

\newsavebox{\buildniceer}
\savebox{\buildniceer}(400,0)[tl]
{
\begin{picture}(400,10)
\put (120,-28){\line (0,1){107}}
\put (140,-28){\line (0,1){107}}
\put (160,-28){\line (0,1){107}}
\put (180,-28){\line (0,1){107}}
\put (200,-28){\line (0,1){107}}

\put (110,20){\line (1,0){102}}
\put (110,30){\line (1,0){102}}
\put (110,40){\line (1,0){102}}
\put (110,50){\line (1,0){102}}
\put (110,60){\line (1,0){102}}
\put (110,70){\line (1,0){102}}
\put (110,10){\line (1,0){102}}
\put (110,0){\line (1,0){102}}
\put (110,-10){\line (1,0){102}}
\put (110,-20){\line (1,0){102}}

\put (117,27){$\bullet$}
\put (117,37){$\bullet$}
\put (117,47){$\bullet$}
\put (117,57){$\bullet$}
\put (117,67){$\bullet$}
\put (137,27){$\bullet$}
\put (137,37){$\bullet$}
\put (137,47){$\bullet$}
\put (137,57){$\bullet$}
\put (157,27){$\bullet$}
\put (157,37){$\bullet$}
\put (157,47){$\bullet$}
\put (177,27){$\bullet$}
\put (177,37){$\bullet$}
\put (197,27){$\bullet$}

\put (87,27){$\scriptstyle L_t$}
\put (87,37){$\scriptstyle L_{t-1}$}
\put (92,50){$\scriptscriptstyle \vdots$}
\put (87,67){$\scriptstyle L_{0}$}
\put (87,17){$\scriptstyle L_{t+1}$}
\put (92,1){$\scriptscriptstyle \vdots$}
\put (87,-11){$\scriptstyle L_{s-t}$}
\put (87,-23){$\scriptstyle M_{t+1}$}

\put (112,-38){$\scriptstyle M_0$}
\put (142,-38){$\scriptstyle \dots$}
\put (167,-38){$\scriptstyle M_{t-1}$}
\put (195,-38){$\scriptstyle M_{t}$}
\end{picture}
}

\begin{equation} \label{pict-in-second-thm}
\hbox{
\begin{picture}(400,60)
\put(40,-5){\usebox{\buildniceer}}
\end{picture}
}
\end{equation}

\vskip .7in

\noindent
Clearly this is a geometric link.  Let
$G_{3,s,t}$ be the Gorenstein scheme obtained by $I_{G_{3,s,t}} = I_{Z_{2,t}}
+ I_{Y_{2,s+1,t+1}}$. One can check geometrically that the components of
$G_{3,s,t}$ are of the form described in the statement of the theorem.  The
simplest way to see this is to use the description of $G_{2,s,t}$
and $Z_{2,t}$.  Since the codimension of $G_{3,s,t}$ is 3 and it is the
intersection of $Z_{2,t}$ and $Y_{2,s+1,t+1}$, each component corresponds to
a {\em pair} of intersection points in (\ref{pict-in-second-thm}), where one
intersection point comes from $Z_{2,t}$ (dots) and one from $Y_{2,s+1,t+1}$
(non-dots), provided these intersection points lie on the same vertical or
horizontal line ($L_i$ or $M_i$) so that the codimension will be 3.  Then it
is a simple matter to verify that the components have the form claimed in
the statement of the theorem.
\end{remark}

\begin{remark} \label{existence-of-basic-configurations} In  Theorems
  \ref{acm-scheme} and \ref{max-gor} we assumed the existence of
  sufficiently general linear forms. This can be guaranteed if, for example,
the field $K$ contains sufficiently many elements  or the polynomial ring
  has dimension $n+1 \geq 2t + 
c$  (in order to construct $Z_{c,t}$) or $n+1 \geq s + c$ (in order to
construct  $G_{c, s, t}$). 

Throughout the rest of the paper it is understood that, whenever one of the
schemes $Z_{c,t}$ or $G_{c, s, t}$  is mentioned, it indeed exists and is a
reduced scheme. 
\end{remark}


\section{A construction of \acm generalized stick figures with arbitrary
Hilbert function} \label{arb-hf-sect} 

In this section we give a construction of an \acm generalized stick figure
having any prescribed Hilbert function in codimension $c$.  The key goal is
to see that they can be viewed as subconfigurations of the generalized stick
figures $Z_{c,t}$ constructed in the previous section, and hence are
automatically  contained in $G_{c,s,t}$. Moreover, we derive a
combinatorial description of the components of the constructed generalized
stick figures.

\begin{definition}
Let $>$ denote the degree-lexicographic order on monomials in the ring $T =
K[z_1,\dots,z_c]$, i.e.\
\[
z_1^{a_1}\cdots z_c^{a_c} > z_1^{b_1}\cdots z_c^{b_c}
\]
 if the first nonzero
coordinate of the vector
\[
\left ( \sum_{i=1}^c (a_i - b_i), a_1 - b_1 ,\dots,a_c - b_c \right )
\]
is positive.  Let $J$ be a monomial ideal.  Let $m_1,m_2$ be monomials in
$T$ of the same degree such that $m_1 > m_2$.  Then $J$ is a {\em lex-segment
ideal} if $m_2 \in J$ implies $m_1 \in J$.
\end{definition}

\begin{notation}
For a graded module $M$ we denote by $a(M)$ the
initial degree:
\[
a(M) := \inf \{ t \in {\mathbb Z} | [M]_t \neq 0 \}.
\]
\end{notation}

\begin{lemma} \label{tool}
Let $J$ be an Artinian lex-segment ideal in $T = K[z_1,\dots,z_c]$ of initial
degree $\alpha$ and for which the $h$-vector ends in degree $t$.  Then there
is a unique decomposition
\[
J = \sum_{j=0}^{\alpha} z_1^j \cdot I_j
\]
where $I_0 \subset I_1 \subset \cdots \subset I_{\alpha-1} \subsetneq
I_{\alpha} = T$,
\begin{equation}\label{shorter-h-vtr}
(z_2,\dots,z_c)^{t+1-j} \subset I_j \subset (z_2,\dots,z_c)^{\alpha-j},
\end{equation}
$$
a(I_j) > \reg (I_{j+1}), 
$$
and $I_j \cap \bar T$ is an Artinian lex-segment ideal in $\bar
T=K[z_2,\dots,z_c]$
$(0 \leq j \leq \alpha-1)$.
\end{lemma}

\begin{proof}
The existence of the decomposition is straightforward.  Note that the
uniqueness comes from the requirement 
that $I_0 \subset I_1 \subset \cdots \subset I_{\alpha-1} \subsetneq
I_{\alpha} = R$. It remains to show the inequality. 

Assume that the regularity of $I_{j+1}$ is $s+1$.  This means
that $z_c^{s+1} \in I_{j+1}$ but $z_c^s \notin I_j$.  Hence
$z_1^{j+1}z_c^{s+1} \in J$ but $z_1^{j+1} z_c^s \notin J$.  Therefore
$z_1^{j} z_2^{s+1} < z_1^{j+1} z_c^s$ is also not contained in $J$, i.e.\
$z_2^{s+1} \notin I_j$. It follows that $a(J_j) > s+1 = \reg (I_j)$ which 
 concludes the proof.
\end{proof}

\begin{remark}
Since the $I_j$ we have produced in Lemma \ref{tool} are Artinian lex-segment
ideals and since $I_j \subset I_k$ whenever $j < k$, their Hilbert functions
will satisfy the hypothesis of Theorem \ref{anyO-seq} below.
\end{remark}

\begin{lemma} \label{hf-lemma}
For any $s\geq 0$, we have
\[
h_{T/J}(s) = \sum_{j=0}^{\alpha-1} h_{\bar T/ I_j \cap \bar T}(s-j)
\]
\end{lemma}

\begin{proof}
It is a straightforward computation.
\end{proof} 

\begin{remark}
The lemma shows that the decomposition of Lemma
\ref{tool} corresponds to the type-vectors of \cite{GHS1} or \cite{GHS2}.
\end{remark}

Lemma \ref{tool} suggests the following definition:

\begin{definition}
Let $\underline{h} = (h_0,h_1,\dots,h_t)$ be an O-sequence.  Choose an integer
$c \geq h_1$ and let $J$ be the Artinian lex-segment ideal in $T =
K[z_1,\dots,z_c]$ such that $\underline{h}$ is the $h$-vector of $T/J$. 
Define $\underline{h}^j$ ($0 \leq j < \alpha$) to be the $h$-vector of
$\bar T / I_j \cap \bar T$, where $\sum_{j=0}^\alpha z_1^j I_j$ is the unique
decomposition according to Lemma \ref{tool}.  We call $(\underline{h}^0,\dots,
\underline{h}^{\alpha-1})$ the {\em decomposition of $\underline{h}$}.
\end{definition}

Note that in \cite{GMR}, \cite{GHS1} or \cite{GHS2} a purely numerical
procedure is given which computes the decomposition of $\underline{h}$, but
does not involve the computation of lex-segment ideals.  However, we need the
relationship to these ideals later on.

Recall that the construction of the scheme $Z_{c,t}$ in Section
\ref{max-hf-sect} involves linear forms $M_i$ and $L_j$ in a polynomial ring
$R$ of dimension $n+1$ over the field $K$.  We are now ready for a
generalization of this construction.

\begin{theorem}\label{anyO-seq}
Let $\underline{h} = (h_0,h_1,\dots,h_v)$ be an O-sequence, where $h_v \neq
0$.  Let $c \geq h_1$, $t \geq v$ be integers.  Suppose
\[
{\mathcal M}_{c,t} = \{ M_0,\dots,M_{t+\lfloor \frac{c-1}{2} \rfloor},
L_0,\dots,L_{t+\lfloor \frac{c-2}{2} \rfloor} \} \subset R
\]
is a set of linear forms such that each subset of of ${\mathcal M}_{c,t}$
consisting of $c+1$ elements is linearly independent (over $K$).  Define the
ideal $I_{c,t} (\underline{h})$ recursively:
\begin{itemize}
\item[] If $c=1$, put $I_{c,t}(\underline{h}):= (M_{t-v}\cdot M_{t-v+1} \cdot
\hdots \cdot M_t)$.

\item[] If $c>1$ put 
\[
I_{c,t}(\underline{h}) :=
\left \{
\begin{array}{ll}
\displaystyle
\bigcap_{j=0}^{\alpha-1}
\left [
I_{c-1,t-j} (\underline{h}^j) + (L_{t-j+\frac{c-2}{2}} )
\right ]
&
\hbox{if $c$ is even} \\
\displaystyle 
\bigcap_{j=0}^{\alpha-1}
\left [
I_{c-1,t-j} (\underline{h}^j) + (M_{t-j+\frac{c-1}{2}} )
\right ] & \hbox{if $c$ is odd}
\end{array}
\right.
\]
where $(\underline{h}^0,\dots,\underline{h}^{\alpha-1})$ is the decomposition
of $\underline{h}$.
\end{itemize}
Then $I_{c,t}(\underline{h})$ defines a reduced \acm subscheme
$Z_{c,t}(\underline{h}) \subset \proj{n}$ of codimension $c$ such that
$Z_{c,t}(\underline{h}) \subset Z_{c,t}$ and the $h$-vector of
$Z_{c,t}(\underline{h})$ is $\underline{h}$.

If each subset of $c+2$ elements of ${\mathcal M}_{c,t}$ is linearly
independent then $Z_{c,t}(\underline{h})$ is a generalized stick figure. 
Moreover, if $\underline{h}' = (h_0',\dots,h_{v'}')$ is an O-sequence where
$h_{v'}' \neq 0$, such that $h_i' \leq h_i$ for all $i$  with $0 \leq i \leq
v' \leq v$, then $Z_{c,t}(\underline{h}') \subset Z_{c,t}(\underline{h})$.
\end{theorem}

\begin{proof}
In view of Theorem \ref{acm-scheme}, it suffices to check that
$Z_{c,t}(\underline{h})$ is arithmetically Cohen-Macaulay,
$Z_{c,t}(\underline{h}) \subset Z_{c,t}$, the $h$-vector of
$Z_{c,t}(\underline{h})$ and the last claim.  We induct on $c \geq 1$.

Let $c=1$.  By comparison with Theorem \ref{acm-scheme} we get
$Z_{c,t}(\underline{h}) \subset Z_{1,t}$ because $I_{Z_{1,t}} =
(M_0\cdot \hdots \cdot M_t)$.

Now let $c>1$.  Using the notation of Lemma \ref{tool}, $\underline{h}^j$ is
the $h$-vector of $\bar T/I_j \cap \bar T$, where $\bar T =
K[z_2,\dots,z_c]$.  Let $\underline{h}^j = (h_0^j,h_1^j,\dots,h_{t_j}^j)$.
Then $h_i^j \leq \binom{c-2+i}{i}$ for all $i \geq 0$.  Moreover, Lemma
\ref{tool} shows that $t_j \leq v-j \leq t-j$.  Since $I_0 \subset I_1
\subset \cdots \subset  I_{\alpha-1}$, the induction hypothesis provides
inclusions of \acm configurations
\[
\begin{array}{ccc}
Z_{c-1,t}(\underline{h}^0) & \subset & Z_{c-1,t} \\
\cup && \cup \\
Z_{c-1,t-1}(\underline{h}^1) & \subset & Z_{c-1,t-1} \\
\cup && \cup \\
\vdots && \vdots \\
\cup && \cup \\
Z_{c-1,t+1-\alpha}(\underline{h}^{\alpha-1}) & \subset & Z_{c-1,t+1-\alpha} \\
&& \cup \\
&& \vdots \\
&& \cup \\
&& Z_{c-1,0}.
\end{array}
\]
Suppose that $c$ is even.  Then the configuration $Z_{c-1,t-j}$ is made up
using only the linear forms 
\[
 M_0,\dots,M_{t-j+ \frac{c-2}{2} },L_0,\dots,L_{t-j+\frac{c-4}{2}}.
\]
Hence the assumption on the set ${\mathcal M}_{c,t}$ and Theorem
\ref{acm-scheme} imply
\[
I_{Z_{c-1,t-j}} : L_{t-j+\frac{c-2}{2}} = I_{Z_{c-1,t-j}}
\]
and hence 
\[
I_{c-1,t-j}(\underline{h}^j) : L_{t-j+\frac{c-2}{2}} =  I_{Z_{c-1,t-j}}.
\]
Therefore Corollary \ref{hyperplane-sects} gives us that the ideal 
\[
(L_{t+\frac{c-2}{2}} \cdot \hdots \cdot L_{t+1-\alpha + \frac{c-2}{2}}) +
\sum_{j=0}^{\alpha-1} I_{c-1,t-j}(\underline{h}^j) \cdot \prod_{i=1}^j
L_{t+1-i+\frac{c-2}{2}}
\]
defines an \acm subscheme having $\underline{h}$ as $h$-vector and is in fact
the ideal $I_{c,t}(\underline{h})$.

If $c$ is odd we conclude similarly because then
\[
I_{c-1,t-j}(\underline{h}^j) : M_{t-j+\frac{c-1}{2}} =
I_{c-1,t-j}(\underline{h}^j).
\]

It remains to show that $I_{c,t}(\underline{h}) \subset
I_{c,t}(\underline{h}')$. Let $J' = \sum_{j=0}^{\alpha'} z_1^j I_j'$ be the
decomposition of the lex-segment ideal $J' \subset T$ with $h$-vector
$\underline{h}'$.  Let $\underline{k}^j$ be the $h$-vector of $\bar T /I_j'
\cap \bar T$.  Lemma \ref{tool} implies that $I_j \subset I_j'$ for all $j$
with $0 \leq j \leq <'$.  Therefore, the induction hypothesis provides
\[
I_{c-1,t-j}(\underline{h}^j) \subset I_{c-1,t-j} (\underline{k}^j)
\ \ \ \hbox{if $0 \leq j < \alpha'$}.
\]
Then using the definition of the ideals we obtain $I_{c,t}(\underline{h})
\subset I_{c,t}(\underline{h}')$ as desired.

Finally, note that if 
\[
\underline{h} = \left ( 1,c,\binom{c+1}{2},\dots,\binom{c-1+t}{t} \right )
\]
is a ``maximal'' $h$-vector then $Z_{c,t}(\underline{h}) = Z_{c,t}$.  Hence
we conclude that all of the configurations $Z_{c,t}(\underline{h})$ that we
obtain for ``smaller'' $h$-vectors are contained in $Z_{c,t}$ as claimed.
\end{proof}

\begin{example}
In Theorem \ref{anyO-seq} we gave the result that if $\underline{h}' =
(h_0',\dots,h_v')$ and $\underline{h} = (h_0,\dots,h_v)$ are O-sequences,
with $v' \leq v \leq t$ and $h_i' \leq h_i$ for all $i$, then
$Z_{c,t}(\underline{h}') \subset Z_{c,t}(\underline{h})$.  We want to stress
here that the value of $t$ must be the same.

Let $\underline{h}' = (1,2)$ and let $\underline{h} = (1,2,2)$.  Then we
obtain
\[
\begin{array}{rcl}
I_{2,1}(\underline{h}') & = & (M_0,L_0) \cap (M_0,L_1) \cap (M_1,L_1) \\
I_{2,2}(\underline{h}') & = & (M_1,L_1) \cap (M_1,L_2) \cap (M_2,L_2) \\
I_{2,2}(\underline{h}) & = & (M_0,L_1) \cap (M_1,L_1) \cap (M_0,L_2) \cap
(M_1,L_2) \cap (M_2,L_2)
\end{array}
\]
In particular, $Z_{2,2}(\underline{h}')$ is contained in
$Z_{2,2}(\underline{h})$ (as expected from Theorem \ref{anyO-seq}) but
$Z_{2,1} (\underline{h}')$ is not contained in $Z_{2,2}(\underline{h})$.
\end{example}

In the following corollary, the use of the term ``extremal Betti numbers''
will be justified in Theorem \ref{big-hul}.

\begin{corollary}\label{max-betti}
The schemes $Z_{c,t}(\underline{h}) \subset \PP^n$ constructed in Theorem
\ref{anyO-seq} have extremal Betti  numbers, i.e.\ their graded Betti numbers
are the same as those of the Artinian lex-segment ideal $J \subset T =
K[z_1,\dots,z_c]$ with $h$-vector $\underline{h}$.
\end{corollary}

\begin{proof} As a preparatory
step we consider the lex-segment ideal  $J \subset T$. The properties of its
decomposition imply 
\[
J : z_1 = I_{0} + \sum_{j=1}^{\alpha} z_1^{j-1} \cdot I_{j} =
\sum_{j=1}^{\alpha} z_1^{j-1} \cdot I_{j}.  
\]
Using Lemmas \ref{hf-lemma} and \ref{tool} we obtain 
\[
\reg (J : z_1) = \reg I_{1} < a(I_{0}). 
\]
Since $(J + z_1 T)/z_1 T = (I_{0} + z_1 T)/z_1 T$, the multiplication by
$z_1$ provides the exact sequence 
\[
0 \to (J : z_1) (-1) \to J \to (I_{0} + z_1 T)/z_1 T \rightarrow 0
\]
where $a((I_{0} + z_1 T)/z_1 T) = a(I_{0}) > \reg (J : z_1)$. We claim
that this 
implies for all $i$ 
\[
\TT^T_i (J, K) \cong \TT^T_i(J : z_1, K) (-1) \oplus \TT^T_i((I_{0} + z_1
T)/z_1 T, K). 
\]
Indeed, the long exact $\TT$-sequence provides maps $\ffi_i: \TT^T_i((I_{0}
+ z_1 T)/z_1 T, K) \to \TT^T_{i-1}(J : z_1, K) (-1)$. If $j < i + a(I_{0})$
then
$[\TT^T_i((I_{0} + z_1 T)/z_1 T, K)]_j = 0$, and if $j \geq i + a(I_{0})
> i + 
\reg (J : z_1)$ then $[\TT^T_{i-1}(J : z_1, K) (-1)]_j = 0$. Therefore $\ffi$
is the zero map and the $\TT$-sequence proves our claim. 

Now we consider $Z_{c,t}(\underline{h})$.  We induct on $c$ and $a(J) =
\alpha$. Again, the case  $c = 1$ is easy. Let $c > 1$. If $\alpha = 1$ then
$J$ contains a linear form which allows us to conclude using induction on the
codimension.  Let $\alpha > 1$. 

Assume that $c$ is even. Let $W$ be defined by 
\[
I_W = I_{c,t}(\underline{h}) :
L_{t +  \frac{c-2}{2}} = \bigcap_{j=1}^{\alpha-1} \left [ I_{c-1,t-j}
(\underline{h}^j) + (L_{t-j+\frac{c-2}{2}} )
\right ].
\]
Then we claim that we have an exact sequence 
\begin{equation}\label{sequence}
0 \to I_{W}(-1) \to I_{c,t}(\underline{h}) \to
\frac{I_{c-1,t}(\underline{h}^0) +  (L_{t +
\frac{c-2}{2}})}{(L_{t+\frac{c-2}{2}})} \to 0. 
\end{equation}
Indeed, the only question is the cokernel.  The sequence clearly holds if we
write
\begin{equation}\label{ideal}
\frac{I_{c,t}(\underline{h}) + (L_{t + \frac{c-2}{2}})
}{(L_{t+\frac{c-2}{2}})}
\end{equation}
 for the cokernel.  But $W$
is \acm by induction, and $Z_{c,t}(\underline{h})$ is \acm by Theorem
\ref{anyO-seq}, so sheafifying and taking cohomology gives that (\ref{ideal})
is a Cohen-Macaulay ideal.  But then the generality of $L_{t +
\frac{c-2}{2}}$ gives that the components of the corresponding scheme are
defined precisely by the ideal given as the cokernel of (\ref{sequence}).

By induction on $\alpha$, $I_W$ has the same Betti numbers as $J : z_1$ and,
by induction on $c$, $I_{c-1,t}(\underline{h}^0)$ has the same Betti numbers
as
$I_{0}$.  Hence we can conclude as above that the Betti numbers of
$I_{c,t}(\underline{h})$ are the sum of the Betti numbers of $I_W(-1)$ and 
\[
\frac{I_{c-1,t}(\underline{h}^0) +  (L_{t +
\frac{c-2}{2}})}{(L_{t+\frac{c-2}{2}})}
\]
 and thus equal to the Betti numbers
of $J$. 

The case where $c$ is odd is handled similarly.
\end{proof}

Some arguments in the proofs of this section could be replaced by results
 from the theory of $k$-configurations (e.g.\ \cite{GHS1}, \cite{GHS2}). But
we prefer to keep the exposition more self-contained.  

We would now like to give an explicit (combinatorial) primary decomposition of
the ideal $I_{c,t}(\underline{h})$.  In order to do this, we will use 
lexicographic order ideals of monomials.  Denote by $S^{(c)}$ the set of
monomials in the polynomial ring $S = K[y_1,\dots,y_c]$.  The {\em reverse
lexicographic order} on $S^{(c)}$ is defined by
\begin{quote}
{\em 
$y_1^{a_1} \cdot \hdots \cdot y_c^{a_c} <_r y_1^{b_1} \cdot \hdots \cdot
y_c^{b_c}$ if the last non-zero coordinate of the vector
\[
(a_1-b_1,\dots,a_c-b_c)
\]
is negative.}
\end{quote}
A non-empty subset ${\mathcal M}$ of $S^{(c)}$ is called an {\em order ideal
of monomials} if $m' \in {\mathcal M}$ and $m|m'$ imply $m \in {\mathcal
M}$.  It is said to be a {\em lexicographic set of monomials} if $m' \in
{\mathcal M}$, $m <_r m'$ and $\deg m = \deg m'$ imply $m \in {\mathcal M}$.
If $\mathcal M$ has both properties, it is called a {\em lexicographic order
ideal of monomials (LOIM)}.  Observe that a LOIM is {\em not} an ideal of $S$.

The {\em lexicographic order} of the set of monomials $T^{(c)}$ in $T =
K[z_1,\dots,z_c]$ is defined by 
\begin{quote}
{\em 
$z_1^{a_1} \cdot \hdots \cdot z_c^{a_c} >_\ell z_1^{b_1} \cdot \hdots \cdot
z_c^{b_c}$ if the first non-zero coordinate of the vector
\[
(a_1-b_1,\dots,a_c-b_c)
\]
is positive.}
\end{quote}
We define a bijective map $\varphi: T^{(c)} \rightarrow S^{(c)}$ by
\[
\varphi(z_1^{a_1}\cdot \hdots \cdot z_c^{a_c}) = y_1^{a_c}\cdot \hdots \cdot
y_c^{a_1}.
\]
Obviously, $\varphi$ preserves the degree.  Moreover, the following is
immediate:

\begin{lemma} \label{rl-lemma}
If $m,m'$ are monomials in $T^{(c)}$ then
\[
m >_\ell m' \ \ \hbox{ if and only if } \ \ \varphi(m) >_r \varphi(m').
\]
\end{lemma}

We also need the following preparatory result.

\begin{lemma} \label{BL-fact}
Let $\underline{h} = (h_0,h_1,\dots)$ be an O-sequence,
where
$c \geq h_1$.  Let ${\mathcal M}_i$ be the smallest (in the ordering $<_r$)
$h_i$ monomials in $S^{(c)}$ of degree $i$.  Then ${\mathcal M} :=
\bigcup_{i=0}^\infty {\mathcal M}_i$ is a LOIM which is called {\em the LOIM
associated to $\underline{h}$} and denoted by $LOIM(\underline{h})$.
\end{lemma}

\begin{proof}
This is a special case of Proposition 1 in \cite{BL2}.  See also
\cite{stanley}.
\end{proof}

Now we can relate the lex-segment ideals in $T$ to the LOIM's in $S$.

\begin{lemma} \label{mainlemma}
Let $J \subset T$ be a lex-segment ideal.  Put $h_i = h_{T/J}(i)$ and
$\underline{h} = (h_0,h_1,\dots)$.   Define
\[
{\mathcal M} := \{ \varphi (m) \ | \ m \in T^{(c)} \backslash J \}.
\]
Then we have ${\mathcal M} = LOIM(\underline{h})$.
\end{lemma}

\begin{proof}
By the definition of $h_i$, the ideal $J$ does not contain exactly
the smallest (with respect to $>_\ell$) $h_i$ monomials in $T$ of degree
$i$.  Therefore, by Lemma \ref{rl-lemma}, $\mathcal M$ contains exactly the
smallest (in the ordering $<_r$) $h_i$ monomials of $S$ of degree $i$.  Since
$\underline{h}$ is an O-sequence, Lemma \ref{BL-fact} shows that $\mathcal M$
is a LOIM.
\end{proof}

\begin{corollary}
Let $\underline{h} = (h_0,\dots,h_t)$ be an O-sequence with the decomposition
\linebreak
$(\underline{h}^0,\dots,\underline{h}^{\alpha-1})$
 and $h_1 \leq c$.  Then we
have 
\[
LOIM (\underline{h}) = \bigcup_{j=0}^{\alpha-1} y_c^j \cdot
LOIM(\underline{h}^j).
\]
(Note that $LOIM(\underline{h}^j) \subset K[y_1,\dots,y_{c-1}]$.)
\end{corollary}

\begin{proof}
Let $J \subset T$ be the lex-segment ideal having $\underline{h}$ as
$h$-vector.  Let $J = \sum_{j=0}^{\alpha-1} z_1^j I_j$ be its decomposition
according to Lemma \ref{tool}.  Since $z_1^\alpha \in J$ we have for a
monomial $z_1^j m$ with $m \in K[z_2,\dots,z_c]$ that 
\[
\begin{array}{rcl}
z_1^j m \notin J & \Leftrightarrow & 0 \leq j < \alpha \hbox{ and } m \notin
I_j \\
& \Leftrightarrow & \varphi(z_1^j m) = y_c^j \varphi(m) \in
LOIM(\underline{h}), \hbox{ where } 0 \leq j < \alpha \hbox{ and } \varphi(m)
\in LOIM(\underline{h}^j) 
\end{array}
\]
by Lemma \ref{mainlemma}.  The claim then follows.
\end{proof}

Now we need some more notation.

\begin{definition}\label{udef}
Let $c \geq 1$, $t \geq 0$ be integers and let $m \in S^{(c)}$ be a monomial
of degree $k \leq t$.  Let $U := K[u_1,\dots,u_{c+2t}]$ be a polynomial
ring.  Write $m$ as
\[
m = y_{e_1} \cdot y_{e_2} \cdot \hdots \cdot y_{e_k} \ \ \hbox{ where } \ \ 
1 \leq e_1 \leq e_2 \leq \dots \leq e_k \leq c.
\]
Define
\begin{itemize}
\item $\displaystyle 
\bar \beta_{c,t} (m)  :=  \{ u_1,\dots,u_{2(t-k)} \} \cup
\bigcup_{j=1}^k \{ u_{e_j+2(j+t-k)-1} , u_{e_j+2(j+t-k)} \} \subset U $

\item $\displaystyle \bar{\mathfrak p}_{c,t} (m)  \hbox{ to be the ideal
generated by }
\{ u_1,\dots,u_{c+2t} \} \backslash \bar \beta_{c,t} (m) $

\item $\displaystyle {\mathfrak p}_{c,t}(m)  \hbox{ to be the ideal
generated by }
\{
\mu (u_i) \ | \ u_i \in \bar {\mathfrak p}_{c,t} (m) \} $ where 
\[ \label{defofmu}
\mu:\{ u_1,\dots,u_{c+2t} \} \rightarrow {\mathcal M}_{c,t} = \{ M_0,\dots,
M_{t + \lfloor \frac{c-1}{2} \rfloor} , L_0,\dots,L_{t+ \lfloor \frac{c-2}{2}
\rfloor} \}
\]
is the bijective map defined by 
\[
\mu (u_i) = 
\left \{
\begin{array}{rl}
\displaystyle M_{\frac{i-1}{2}} & \hbox{ if $i$ is odd}; \\
\displaystyle L_{\frac{i-2}{2}} & \hbox{ if $i$ is even.}
\end{array}
\right.
\]
\end{itemize}
\end{definition}

\begin{theorem} \label{decompofh}
Let $\underline{h} = (h_0,\dots,h_t)$ be an O-sequence, where $h_1 \leq c$. 
Then we have 
\[
I_{c,t} (\underline{h}) = \bigcap_{m \in LOIM (\underline{h})} {\mathfrak
p}_{c,t} (m).
\]
\end{theorem}

\begin{proof}
We induct on $c \geq 1$.  First let $c=1$.  Then there is an integer $v$ such
that $0 \leq v \leq t$ and 
\[
h_i = \left \{
\begin{array}{rl}
1 & \hbox{ if $0 \leq i \leq v$} \\
0 & \hbox{ if $v < i \leq t$}.
\end{array}
\right.
\]
Thus we have 
\[
I_{c,t}(\underline{h}) = (M_{t-v} \cdot \hdots \cdot M_t ) \ \ \hbox{ and } \
\ LOIM(\underline{h}) = \{ 1,y_1, \dots, y_1^v \}.
\]
For $0 \leq k \leq v$ we obtain 
\[
\bar \beta_{c,t} (y_1^k) = \{ u_1,\dots,u_{2(t-k)} \}
\cup
\{ u_{2(t-k)+2} ,\dots, u_{2t+1} \},
\]
and thus $\bar {\mathfrak p}_{c,t} (y_1^k) =  (u_{2(t-k)+1} )$, and
${\mathfrak p}_{c,t} (y_1^k) = M_{t-k}$.  It follows that
\[
\bigcap_{m \in LOIM(\underline{h})} {\mathfrak p}_{c,t}(m) =
\bigcap_{k=0}^v (M_{t-k}) = I_{c,t}(\underline{h})
\]
as claimed.

Now let $c>1$.  Let $(\underline{h}^0,\dots,\underline{h}^{\alpha-1})$ be the
decomposition of $\underline{h}$.  Assume that $c$ is even.  Then we know by
Theorem \ref{anyO-seq} and the induction hypothesis that
\[
\begin{array}{rcl}
\displaystyle I_{c,t}(\underline{h}) & = & \displaystyle 
\bigcap_{j=0}^{\alpha-1} \left [ I_{c-1,t-j} (\underline{h}^j) +
(L_{t-j+\frac{c-2}{2}}) \right ] \\
& = &
\displaystyle 
\bigcap_{j=0}^{\alpha-1} \left [ \bigcap_{m \in LOIM(\underline{h}^j)}
{\mathfrak p}_{c-1,t-j}(m) + (L_{t-j+\frac{c-2}{2}}) \right ].
\end{array}
\]
Since $LOIM(\underline{h}) = \bigcup_{j=0}^{\alpha-1} y_c^j \cdot LOIM
(\underline{h}^j)$, we are done if we can show

\noindent \underline{Claim}: Let $0 \leq j < \alpha$ and let $m \in
LOIM(\underline{h}^j)$.  Then 
\[
{\mathfrak p}_{c,t} (y_c^j \cdot m) = {\mathfrak p}_{c-1,t-j}(m) +
(L_{t-j+\frac{c-2}{2}}).
\]

By Lemma \ref{tool} we know that $k := \deg m \leq t-j$.    Thus we get
\[
\bar \beta_{c,t} (y_c^j \cdot m) = \bar \beta_{c-1,t-j}(m) \cup
\{ u_{c+2(t-j)+1},\dots,u_{c+2t} \}.
\]
Because $\bar {\mathfrak p}_{c-1,t-j}(m)$ is generated by $\{
u_1,\dots,u_{c-1+2(t-j)} \} \backslash \bar \beta_{c-1,t-j}(m)$ and $\bar
{\mathfrak p}_{c,t} (y_c^j \cdot m)$ is generated by $\{ u_1,\dots, u_{c+2t}
\}
\backslash \bar \beta_{c,t}(y_c^j \cdot m)$, it follows that 
$\bar {\mathfrak p}_{c,t} (y_c \cdot m) = \bar {\mathfrak p}_{c-1,t-j}(m) +
(u_{c+2(t-j)})$, and so
\[
{\mathfrak p}_{c,t} (y_c^j \cdot m) = {\mathfrak p}_{c-1,t-j} (m) +
(L_{t-j+\frac{c-2}{2}} ),
\]
proving the claim.

If $c$ is odd we conclude similarly.
\end{proof}


\section{Gorenstein configurations for any SI-sequence}
\label{tie-together}

In this section we prove one of the main results of the paper, namely that
for every SI-sequence $\underline{h}$ there exists a reduced, arithmetically
Gorenstein union of linear varieties whose $h$-vector is $\underline{h}$. 
Furthermore, we characterize the linear varieties which are components of
this arithmetically Gorenstein configuration, essentially giving the primary
decomposition.

We begin by introducing a notation for a shift among the variables $u_i$:

\begin{definition}
We define the injective map $\tau: \{ u_1,u_2,\dots \} \rightarrow
\{ u_1,u_2,\dots  \}$ by $\tau(u_i) = u_{i+1}$ for all $i \geq 1$.
Abusing notation, we denote by $\tau({\mathfrak p}_{c-1,t} (m))$ the ideal
generated by
\[
\{ \mu(\tau(u_i)) \ | \ u_i \in \bar {\mathfrak{p}}_{c-1,t}(m) \}
\]
and by $\tau(Z_{c-1,t})$ the scheme defined by the intersection of the ideals
$\tau (\mathfrak{p})$ where $\mathfrak p$ is a minimal prime ideal of
$I_{c-1,t}$ (using Theorems \ref{anyO-seq} and \ref{decompofh}).
\end{definition}

In the next technical lemma we collect some facts that we will need later on.

\begin{lemma} \label{fax-lemma8}
Let $c \geq2$, $t \geq 1$ be integers.  Then we have

\begin{itemize}

\item[(a)] $\tau (Z_{c-1,t-1}) \subset G_{c-1,s,t}$ for all $s \geq 2t$.

\item[(b)] If $m = y_{e_1}\cdot \hdots \cdot y_{e_k}$ where $1 \leq e_1 \leq
\dots \leq e_k \leq c-1$ and $k \leq t$ then

\begin{itemize}
\item[(i)] $u_1 \in \bar {\beta}_{c-1,t} (m)$ if and only if $k < t$.

\item[(ii)] $u_{c-1+2t} \in \bar {\beta}_{c-1,t} (m)$ if and only if
$e_k = c-1$.

\item[(iii)] If $y_{c-1}$ divides $m$ then $\displaystyle
\mathfrak{p}_{c-1,t}(m) = \mathfrak{p}_{c-1,t-1} \left ( \frac{m}{y_{c-1}}
\right )$.

\end{itemize}
\end{itemize}
\end{lemma}

\begin{proof}
(a) Comparing Theorems \ref{acm-scheme} and \ref{max-gor} we observe that
$\tau (Z_{c-1,t-1})$ is just defined by $B_{c-1,s,t-1}$ and $B'_{c-1,s,t-1}$
respectively, depending on the divisibility of $c-1$ by $2$.

(b) The first two claims follow immediately from the definition of $\bar
\beta_{c-1,t}(m)$  For (iii) it suffices to note that $e_k = c-1$
implies $\{ u_{c-2+2t},u_{c-1+2t} \} \subset \bar \beta_{c-1,t}(m)$, and thus
\[
\{ u_1,\dots,u_{c-1+2t} \} \backslash \bar \beta_{c-1,t} (m) = \{
u_1,\dots,u_{c-1+2(t-1)} \} \backslash \bar \beta_{c-1,t-1} \left (
\frac{m}{y_{c-1}} \right )
\]
and the assertion follows.
\end{proof}

Now we are ready for the announced construction of the arithmetically
Gorenstein schemes and a description of their irreducible components.

\begin{theorem}\label{maingorthm}
Let $\underline{h} = (h_0,h_1,\dots,h_s)$ be an SI-sequence, where $h_s =
1$.  Let $c \geq \max \{h_1, 2\}$ be an integer and let $t = \min \{ i |
h_i \geq h_{i+1} 
\}$.  Put $\underline{g} = (g_0,\dots,g_t)$, where $g_i = h_i - h_{i-1}$. 
Suppose 
\[
{\mathcal N}_{c-1,s+1,t} = \{ M_0,\dots,M_{t+\lfloor \frac{c-2}{2} \rfloor }
,L_0,\dots,L_{s-t+\lfloor \frac{c-1}{2} \rfloor } \} \subset R
\]
is a subset of linear forms such that each subset of $c+1$ elements is
linearly independent.  Define the ideal $J_c(\underline{h}) :=
I_{c-1,t}(\underline{g}) + I_Y$ where $Y = G_{c-1,s+1,t} \backslash
Z_{c-1,t}(\underline{g})$.
  Then
$J_c(\underline{h})$ is a Gorenstein ideal in $R$ defining a reduced,
arithmetically Gorenstein subscheme
$G_c(\underline{h}) \subset \proj{n}$ of codimension $c$ having
$\underline{h}$ as $h$-vector.  Furthermore, $\mathfrak p$ is a minimal prime
ideal of $J_c(\underline{h})$ if and only if it is generated by a subset of
${\mathcal N}_{c-1,s+1,t}$ consisting of $c$ elements and contains exactly
one minimal prime ideal of $I_{c-1,t}(\underline{g})$.
\end{theorem}

\begin{proof}
Notice that $Z_{c-1,t}$ is always contained inside $G_{c-1,s+1,t}$, thanks
to Theorem \ref{acm-scheme} and Theorem \ref{max-gor}.  However, even if
$\underline{h}$ is maximal and so $Z_{c-1,t} = Z_{c-1,t}(\underline{h})$
(cf.\ Theorem \ref{anyO-seq}), it is not true (if $c$ is odd) that
$G_{c,s,t} = G_c(\underline{h})$, or even that 
$Z_{c,t} \subset G_c(\underline{h})$.  See also Remark \ref{compare-constr}.

(Step I) According to Theorem \ref{anyO-seq}, we know that
$Z_{c-1,t}(\underline{g}) \subset Z_{c-1,t}$.  By Theorem \ref{max-gor}, 
$G_{c-1,s+1,t}$ is a generalized stick figure 
containing $Z_{c-1,t}$.  Thus $J_c(\underline{h})$ is the sum of the
geometrically linked ideals $I_{c-1,t}(\underline{g})$ and $I_Y$.  The fact
that $G_c(\underline{h})$ is reduced comes from the fact that $G_{c-1,s+1,t}$
is a generalized stick figure, using Remark \ref{can-add}.

We have
\[
\underline{h} = (1,c, h_2, \dots,h_{t-1}, h_t,
h_t, \dots, h_t,h_t,h_{t+1},\dots,h_{s-2}, c,1).
\]
Note that the 1's occur in degrees 0 and $s$ and the last $h_t$ occurs in
degree $s-t$.  Note also that $s \geq 2t$.
Consider the sequences $\underline{d} = (d_i)$, $\underline{g} = (g_i)$,
$\underline{c} = (c_i)$, $\underline{g}' = (g'_i)$, where
\[
\begin{array}{rcll}
d_i & = & h_i - h_{i-1} & \hbox{for $0 \leq i \leq s+1$  (some of the
entries are 0 or}\\
&&& \hbox{negative and the sequence is anti-symmetric);}\\
g_i & = & 
\left \{
\begin{array}{l}
h_i - h_{i-1}  \\
0
\end{array}
\right.
& 
\! \! \!
\begin{array}{l}
\hbox{for $0 \leq i \leq t$  (this is an
O-sequence);}\\
\hbox{for $i>t$;}
\end{array}
\\
c_i & = &
\left \{
\begin{array}{l}
\binom{i+c-2}{c-2}  \\
\binom{t+c-2}{c-2}  \\
\binom{s-i+c-1}{c-2} \\
\end{array}
\right.
&
\! \! \!
\begin{array}{l}
\hbox{for $0 \leq i \leq t$}; \\
\hbox{for $t \leq i \leq s-t+1$};\\
\hbox{for $s-t+1 \leq i \leq s+1$}
\end{array}
\\
&&& \hbox{(note that $c_i =  c_{s+1-i}$ for $0 \leq i \leq s+1$);} \\
g'_i & = & c_{s+1-i} - g_{s+1-i} & \hbox{for $i \geq 0$}.
\end{array}
\]

Note that $\underline{c}$ is the $h$-vector of $G_{c-1,s+1,t}$.  By Lemma
\ref{comp-tool} we know that $Y$ has $h$-vector $\underline{g}'$, and again by
Lemma  \ref{comp-tool} we know that the first difference of the $h$-vector of
$G_c(\underline{h})$ is given by
\[
g_i + g'_i - c_i.
\]

We now claim that $g_i + g'_i - c_i = d_i$.
\begin{itemize}
\item For $0 \leq i \leq t$ we have $g_i = d_i$ and $g'_i = c_i$.

\item For $t+1 \leq i \leq s-t$ we have $c_i = \binom{t+c-2}{c-2}$, $g_i = 0 =
g_{s+1-i}$ and $g'_i = c_{s+1-i} - 0 = c_i$.  However, in this range we have
$d_i = 0$ by the symmetry of $\underline{h}$.

\item For $s-t+1 \leq i \leq s+1$ we have $g_i = 0$, $g'_i = c_{s+1-i} -
g_{s+1-i} = c_i - g_{s+1-i}$.  Hence
\[
\begin{array}{rcl}
g_i + g'_i - c_i & = & -g_{s+1-i} \\
& = & - d_{s+1-i} \\
& = & d_i
\end{array}
\]
\end{itemize}

\noindent It follows that $G_c(\underline{h})$ has the desired $h$-vector
$\underline{h}$.  

\medskip

(Step II) 
It remains to establish the claim on the minimal prime ideals of
$J_c(\underline{h})$.  Since $G_{c-1,s+1,t}$ is a generalized stick figure
and $J_c(\underline{h})$ is unmixed of codimension $c$, an ideal $\mathfrak
p$ is a minimal prime of $J_c(\underline{h})$ if and only if $\codim
(\mathfrak{p}) = c$ and $\mathfrak{p} = \mathfrak{p}_1 + \mathfrak{p}_2$ for
some minimal prime ideal $\mathfrak{p}_1$ of $I_{c-1,t}(\underline{g})$ and
some minimal prime ideal $\mathfrak{p}_2$ of $I_Y$.  Thus the assumption on
${\mathcal N}_{c-1,s+1,t}$ implies that each minimal prime $\mathfrak p$ of
$J_c(\underline{h})$ is generated by a subset of ${\mathcal N}_{c-1,s+1,t}$
consisting of $c$ elements.  

Suppose that $\mathfrak p$ contains another minimal prime $\mathfrak{p}_1'
\neq \mathfrak{p}_1$ of $I_{c-1,t}(\underline{g})$.  Then $\mathfrak{p}_1,
\mathfrak{p}_1', \mathfrak{p}_2$ are pairwise distinct minimal prime ideals
of $G_{c-1,s+1,t}$ such that $\mathfrak{p}_1 + \mathfrak{p}_1' +
\mathfrak{p}_2 = \mathfrak{p}$ has codimension $c$.  This contradicts the
fact that $G_{c-1,s+1,t}$ is a generalized stick figure.  Thus we have shown
the necessity of the conditions for being a minimal prime.

\medskip

(Step III)
In order to prove  sufficiency, let $\mathfrak{p}_1 =
\mathfrak{p}_{c-1,t}(m)$ be a minimal prime of $I_{c-1,t}(\underline{g})$,
where $m \in LOIM(\underline{g})$.  We must show that for $H \in {\mathcal
N}_{c-1,s+1,t} \backslash \mathfrak{p}_1$, if $\mathfrak{p}_1$ is the only
minimal prime of $I_{c-1,t}(\underline{g})$ which is contained in
$\mathfrak{p}_1 + (H)$ then $\mathfrak{p}_1 + (H)$ is a minimal prime of
$J_c(\underline{h})$.

Our proof will make use of the following observation.  Given $H$ and
$\mathfrak{p}_1$ as in the last paragraph, consider a specific minimal prime
ideal
$\mathfrak{q}$ of
$I_{G_{c-1,s+1,t}}$ such that $\mathfrak{p}_1 +(H) = \mathfrak{p}_1 +
\mathfrak{q}$.  Then we know by Step II that $\mathfrak{p}_1 + (H)$ is a
minimal prime ideal of $J_c(\underline{h})$ if and only if $\mathfrak{q}$ is
not a minimal prime ideal of $I_{c-1,t}(\underline{g})$ (i.e.\ $\mathfrak{q}$
is a minimal prime of $I_Y$).

We will distinguish several cases.

\medskip

\underline{Case 1.}
Suppose that 
$H \in {\mathcal N}_{c-1,s+1,t} \backslash {\mathcal M}_{c-1,t}$ (see the
statement of Theorem \ref{anyO-seq}), i.e.\
\[
H \in \left \{ L_{t+1+\lfloor \frac{c-3}{2} \rfloor},\dots,L_{s+1-t+\lfloor
\frac{c-3}{2} \rfloor} \right \}.
\]
We already know that $\mathfrak{p}_1$ is a minimal prime ideal of
$I_{G_{c-1,s+1,t}}$.  Write $\mathfrak{p}_1$ as in Theorem \ref{max-gor} and
denote by $\mathfrak{q}'$ the ideal generated by the first $c-2$ linear
generators of $\mathfrak{p}_1$.  Then $\mathfrak{q} = \mathfrak{q}' + (H)$ is
a minimal prime ideal of $I_{G_{c-1,s+1,t}}$ -- in fact, it is a primary
component of $C_{c-1,s+1,t}$ (if $c-1$ is even) or $C'_{c-1,s+1,t}$ (if $c-1$
is odd).  Since $H \notin {\mathcal M}_{c-1,t}$, $\mathfrak{q}$ must be a
minimal prime of $I_Y$.  Thus $\mathfrak{p}_1 + (H) = \mathfrak{p}_1 +
\mathfrak{q}$ is a primary component of $J_c(\underline{h})$.

\medskip

\underline{Case 2.}
Suppose $H \in {\mathcal M}_{c-1,t}$, i.e.\ $H = \mu(u)$ for some $u \in \{
u_1,\dots,u_{c-1+2t} \}$.  As preparation we need a description of the subsets
of ${\mathcal U} = \{ u_1,\dots,u_{c-1+2t} \}$ having the form $\bar
\beta_{c-1,t}(m)$ for some monomial $m \in S^{(c-1)}$ (see Definition
\ref{udef}).  To this end we call a subset $W$ of $\mathcal U$ {\em
consecutive} if $W = \{ u_i,u_{i+1},\dots,u_j \}$ for some $1 \leq i \leq j
\leq c-1+2t$.  If $i>1$ then $u_{i-1}$ is called the {\em predecessor} of $W$,
and if $j < c-1+2t$ then $u_{j+1}$ is called the {\em successor} of $W$.

Now let $W$ be an arbitrary non-empty subset of $\mathcal U$.  It is clear
that $W$ has a unique decomposition $W = W_1 \cup \dots \cup W_p$ where
$W_1,\dots,W_p$ are consecutive subsets, the successor of $W_i$ is not
contained in $W$ if $1 \leq i < p$, the predecessor of $W_i$ is not contained
in $W$ if $1 < i \leq p$, and every element of $W_i$ is smaller (in the order
$<_r$) than every element of $W_{i+1}$ if $1 \leq i < p$.  Using this
notation, the definition of $\bar \beta_{c-1,t}(m)$ (Definition \ref{udef})
implies:

\medskip

\noindent \underline{Observation}:  {\em There is a monomial $m \in S^{(c-1)}$
such that $W = \bar \beta_{c-1,t} (m)$ if and only if
$|W| = 2t$ and the cardinality $|W_i|$ is even for all $i$ with $1 \leq i \leq
p$. }

\medskip

Returning to our previous notation, let $m \in S^{(c-1)}$ again denote a
monomial such that $\mathfrak{p}_1 = \mathfrak{p}_{c-1,t}(m)$.  As above,
write $W = W_1 \cup \dots \cup W_p$ for the unique decomposition of $W :=
\bar \beta_{c-1,t}(m)$ into consecutive subsets.  Since $H = \mu(u) \notin
\mathfrak{p}_1$, we obtain $u \in W = \bar \beta_{c-1,t}(m)$, i.e.\ $u \in
W_\ell$ for some $\ell$ ($1 \leq \ell \leq p$).  Write $W_\ell \backslash \{ u
\} = \tilde W \cup \bar W$ where $\tilde W$ and $\bar W$ are both consecutive
unless one of them is empty, and each element of $\tilde W$ is smaller than
each element of $\bar W$.

We have to distinguish further cases:

\medskip

\underline{Case 2.1}: Suppose that $|\tilde W|$ is odd and that $u_1 \notin
W$ if $\ell=1$.  Then $W_\ell$ has a predecessor, say $\bar u$.  Since $|\bar
W|$ is even, the observation above implies the existence of a monomial $m' \in
S^{(c-1)}$ of degree $\leq t$ such that 
\[
(W \backslash \{ u \}) \cup \{ \bar u \} = \bar \beta_{c-1,t}(m').
\]

The prime ideal $\mathfrak{q} := \mathfrak{p}_{c-1,t}(m')$ is a minimal prime
of $I_{Z_{c-1,t}}$, and thus of $I_{G_{c-1,s+1,t}}$.  Moreover, we have by
construction
\[
\mathfrak{p}_1 + \mathfrak{q} = \mathfrak{p}_1 + (\mu(u)) = \mathfrak{p}_1 +
(H).
\]
Hence we are done, as explained at the beginning of Step (III).

\medskip

\underline{Case 2.2}:
Suppose that $|\tilde W|$ is even and that $u_{c-1+2t} \notin W$ if
$\ell=p$.  Then $W_\ell$ has a successor, say $\bar u$.  Since $|\tilde W|$
is even, we can find as above a monomial $m' \in S^{(c-1)}$ such that 
\[
(W \backslash \{ u \}) \cup \{ \bar u \} = \bar \beta_{c-1,t}(m')
\]
and $\mathfrak{q} = \mathfrak{p}_{c-1,t}(m')$ is a minimal prime of
$I_{G_{c-1,s+1,t}}$.  We conclude as in Case 2.1.

\medskip

\underline{Case 2.3}:
Suppose that $|\tilde W|$ is odd, $u_1 \in W$ and $\ell=1$.  Then $u_1 \in W
\backslash \{ u \}$ (since  if $u=u_1$ then $|\tilde W|$ is not odd) 
and $p \geq 2$.  

\underline{Case 2.3.1}
Suppose that $u_{c-1+2t} \notin W$.
Then we form the subset $X \subset {\mathcal U}$ by replacing $u_1$ by
$u_{c-1+2t}$ and all $u_i \in (W \backslash \{ u \}) \cup \{ u_{c-1+2t} \}$
with $i>1$ by $u_{i-1}$.  The observation above shows that there is a
monomial $m' \in S^{(c-1)}$ of degree $\leq t$ such that $X = \bar
\beta_{c-1,t}(m')$.  Since $u_{c-1+2t} \in X$, we obtain by Lemma
\ref{fax-lemma8} that $y_{c-1}$ divides $m'$.  Using Lemma \ref{fax-lemma8}
again, we conclude that 
\[
\mathfrak{q} := \mathfrak{p}_{c-1,t}(m') =
\mathfrak{p}_{c-1,t-1} \left ( \frac{m'}{y_{c-1}} \right )
\]
is a minimal prime ideal of $I_{Z_{c-1,t-1}}$.  Hence Lemma \ref{fax-lemma8}
shows that $\tau(\mathfrak{q})$ is a minimal prime of $I_{G_{c-1,s+1,t}}$. 
It follows that $W \cap \tau(X) = W \backslash \{ u,u_1 \}$, and thus
\[
\mathfrak{p}_1 + \tau(\mathfrak{q}) = \mathfrak{p}_1 + (H)
\]
is a minimal prime ideal of $J_c(\underline{h})$.

\medskip

\underline{Case 2.3.2}: Suppose that $u_{c-1+2t} \in W$.
Let $\bar u$ be the predecessor of $W_p$.  As in Case 2.3.1, construct $X$
out of $(W \backslash \{ u \}) \cup \{ \bar u \}$ by replacing all $u_i \in
(W \backslash \{u\}) \cup \{ \bar u \}$ by $u_{i-1}$ if $i>1$ and by
$u_{c-1+2t}$ if $i=1$.  Again we get the existence of a monomial $m' \in
S^{(c-1)}$ such that $X = \bar \beta_{c-1,t}(m')$ and 
\[
\mathfrak{q} = \mathfrak{p}_{c-1,t}(m') = \mathfrak{p}_{c-1,t-1} \left (
\frac{m'}{y_{c-1}} \right )
\]
is a minimal prime of $I_{Z_{c-1,t-1}}$, and thus $\mathfrak{p}_1 +
\tau(\mathfrak{q}) = \mathfrak{p}_1 +(H)$ is a minimal prime of
$J_c(\underline{h})$.

\medskip

\underline{Case 2.4}: Suppose that $|\tilde{W}|$ is even, $u_{c-1+2t} \in W$
and $\ell=p$.
Then $u_{c-1+2t} \in \bar W$ and $|\bar W|$ is odd.  (Note that 
$u = u_{c-1+2t}$ is impossible since otherwise $\tilde W = W_\ell \setminus
\{ u_{c-1+2t} \}$ has odd cardinality.)

\medskip 

\underline{Case 2.4.1}:
Suppose that $u_1 \notin W$.
Similarly as above we conclude that there is a monomial $m' \in S^{(c-1)}$
such that $\bar \beta_{c-1,t}(m') = (W \backslash \{ u \} ) \cup \{ u_1 \}$,
and
\[
 \mathfrak{q} = \mathfrak{p}_{c-1,t}(m') =
\mathfrak{p}_{c-1,t-1}\left ( \frac{m'}{y_{c-1}} \right )
\]
is a minimal prime of $I_{Z_{c-1,t-1}}$, and thus $\mathfrak{p}_1 +
\tau(\mathfrak{q}) = \mathfrak{p}_1 + (H)$ is a minimal prime of
$J_c(\underline{h})$.

\medskip

\underline{Case 2.4.2}:
Suppose that $u_1 \in W$.
Let $\bar u$ be the successor of $W_1$.  Then there is a monomial $m' \in
S^{(c-1)}$ such that $y_{c-1}$ divides $m'$, $\bar \beta_{c-1,t}(m') = (W
\backslash \{ u \}) \cup \{ \bar u \}$ and 
\[
\mathfrak{p}_1 + \tau \left ( \mathfrak{p}_{c-1,t-1} \left (
\frac{m'}{y_{c-1}} \right ) \right ) = \mathfrak{p}_1 + (H)
\]
is a minimal prime of $J_c(\underline{h})$.

Therefore we have shown the claim in all cases, and the proof of Theorem
\ref{maingorthm} is complete.
\end{proof}

\begin{remark}\label{compare-constr} We conjecture that the \aG
configurations produced in Theorem
\ref{maingorthm} are in fact themselves generalized stick figures.  The
obstruction to proving this comes from the possibility of situations like the
hypothetical one described in the proof of Theorem \ref{max-gor} (see page
\pageref{obst-to-gen-sf}). There we were able to conclude because of the very
explicit nature of the primary decomposition, which is missing here.  Note
that in codimension 3 the result is known \cite{GM5}.

This conjecture would follow easily if we knew, for example, that
$G_c(\underline{h}) \subset G_{c,s,t}$.  There is a subtle point, to begin
with.  The scheme $G_c(\underline{h})$ is constructed using the set 
\[
{\mathcal
N}_{c-1,s+1,t} = \{ M_0,\dots,M_{t+\lfloor \frac{c-2}{2}
\rfloor},L_0,\dots,L_{s+1-t+\lfloor \frac{c-3}{2} \rfloor} \}
\]
 whereas
$G_{c,s,t}$ is constructed using 
\[
{\mathcal N}_{c,s,t} = \{
M_0,\dots,M_{t+\lfloor \frac{c-1}{2} \rfloor},L_0,\dots,L_{s-t+\lfloor
\frac{c-2}{2}\rfloor} \}.
\]
These agree if and only if $c$ is even.  In order to remedy this situation
we can proceed as we did in Theorem \ref{max-gor} and define
\[
G'_{c,s,t} = 
\left \{
\begin{array}{l}
G_{c,s,t} \hbox{ if $c$ is even};\\ 
\hbox{The scheme defined by the ideal which is obtained from $I_{G_{c,s,t}}$
}\\
\hbox{by re-naming $M_{t+\frac{c-1}{2}}$ as $L_{s-t+\frac{c-1}{2}}$.}
\end{array}
\right.
\]
With this re-naming, ${\mathcal N}_{c,s,t}$ becomes ${\mathcal
N}_{c-1,s+1,t}$, and of course $G'_{c,s,t}$ is still arithmetically
Gorenstein.

The ``right'' question is thus whether $G_c(\underline{h}) \subset
G'_{c,s,t}$.  In fact, if $\underline{h}$ is maximal, i.e.\ if 
\[
h_i = \left \{
\begin{array}{rl}
\binom{c-1+i}{c-1} & \hbox{if $0 \leq i \leq t$}; \\
\binom{c-1+t}{c-1} & \hbox{if $t \leq i \leq s-t$};\\
\binom{s-i+c-1}{c-1} & \hbox{if $s-t \leq i \leq s$}
\end{array}
\right.
\]
then it follows from the constructions of Theorem \ref{max-gor} and Theorem
\ref{maingorthm} that $G_c(\underline{h}) = G'_{c,s,t}$.

However, we can not expect, in general, that $G_c(\underline{h}) \subset
G'_{c,s,t}$.  The reason is the following. 
Suppose we have $Z_{c-1,t}(\underline{g}) \subsetneq Z_{c-1,t}
\subset G$, where $G$ is the suitable \aG scheme from the construction.  The
components of $G_{c,s,t}$ are obtained by intersecting a component of
$Z_{c-1,t}$ with a component of $G \backslash Z_{c-1,t}$.  The
components of $G_c(\underline{h})$ are obtained by intersecting a component of
$Z_{c-1,t}(\underline{g})$ with a component of $G \backslash
Z_{c-1,t}(\underline{g})$.  Hence there are components of
$G_c(\underline{h})$ which consist of the intersection of two components of
$Z_{c-1,t}$, and no such intersection is a component of $G'_{c,s,t}$.

For a specific counterexample, let $\underline{h} = (1,3,5,3,1)$, so $s=4$,
$t=2$ and we put $c=3$.  then $(L_0,M_0,L_1)$ is a minimal prime ideal of
$J_3(\underline{h})$ but not of $I_{G'_{3,4,2}}$.  Hence $G_3(\underline{h})
\not \subset G'_{3,4,2}$.
\end{remark}


\section{The Subspace Property} \label{subspace-sect}

In this section we want to show that the \aG schemes constructed in the
previous section have the Weak Lefschetz Property.  Passing to the Artinian
reduction we would lose our useful combinatorial description.  Thus, we
introduce the so-called {\em subspace property}.  It is even defined for
\acm schemes and does not require taking consecutive hyperplane sections. 
However, for \aG subschemes the subspace property implies the Weak Lefschetz
Property.  Thus, we conclude by showing that our \aG schemes do possess the
subspace property.

Recall (cf.\ Remark \ref{def-reg} ) that $\reg(X)$ denotes the
Castelnuovo-Mumford regularity of the scheme $X$.

\begin{notation}
For a graded module $M$ we denote by $a(M)$ the
initial degree:
\[
a(M) := \inf \{ t \in {\mathbb Z} | [M]_t \neq 0 \}.
\]
\end{notation}

\begin{definition}\label{subspaceprop}
A homogeneous ideal $I \subset R = K[x_0,\dots,x_n]$ of codimension $c$ is
said to have the {\em subspace property} if there is a linear form $\ell \in
R$ such that $(I+\ell R)/ \ell R$ is a perfect ideal in $R/ \ell R$ of
codimension $c-1$ and the initial degree of $(I:_R \ell)/ I$ is at least
$\frac{\reg (I) -1}{2}$.  $X \subset \proj{n}$ is said to have the subspace
property if $I_X$ does.
\end{definition}

\begin{remark}
The condition that $(I+\ell R)/ \ell R$ have codimension $c-1$ in $R/ \ell R$
means that $\ell$ has to vanish on some components of the scheme defined by
$I$.  The conditions that $(I+\ell R)/ \ell R$ is perfect and that the initial
degree of $(I:_R \ell)/ I$ is ``not too small'' are what is difficult to
verify, in general.

The main example for us is the following.  Let $I_S$ be
the saturated ideal of an \acm subscheme $S$ of codimension $c-1$ in
$\proj{n}$.  Let $I_C$ be the saturated ideal of an \acm subscheme $C$ of $S$
of codimension $c$ in $\proj{n}$. Let $\ell$ be a linear form not vanishing on
any component of $S$ or of $C$.  By Lemma \ref{KMMNP-lemma}, $I = I_S + \ell
\cdot I_C$ is the saturated ideal of an \acm subscheme of codimension
$c$ in $\proj{n}$, and this subscheme consists of the union of $C$ and a
hyperplane section of $S$.  The ideal $I + \ell R/\ell R$ is just the ideal of
this hyperplane section {\em inside the hyperplane}, and hence is perfect and
has codimension $c-1$ inside the hyperplane as required.  The ideal $I:_R \ell
$ is just $I_C$, so if we start with $C$ which is rather large inside $S$,
the condition on the regularity is easy to verify.
\end{remark}

\begin{lemma}\label{artin-sub-iff-wlp}
Let $R/I$ be an Artinian Gorenstein algebra.  Then $R/I$ has the weak
Lefschetz property with respect to a linear form $g$ (cf.\ Definition
\ref{WSP}) if and only if  $I$ has the subspace property with respect to $g$.
\end{lemma}

\begin{proof}
Assume that $R/I$ has the weak Lefschetz property with respect to $g$.  The
fact that we get a perfect ideal comes from the Artinian property.   Let
\[
(h_0,h_1,\dots,h_{t-1},h_t,\dots,h_t, \dots, h_1, h_0)
\]
 be the
$h$-vector of $R/I$ where $h_{t-1} < h_t$ and let $s := \reg (R/I)$ (so the
second $h_0 =1$ occurs in degree $s$).  Consider the exact sequence
\begin{equation} \label{exseq}
0 \rightarrow \frac{I:g}{I} (-1) \rightarrow \frac{R}{I}(-1)
\stackrel{g}{\longrightarrow} \frac{R}{I} \rightarrow \frac{R}{I+gR}
\rightarrow 0.
\end{equation}
Since the multiplication by $g$ on $R/I$ is injective if $h_{i-1} \leq h_i$,
we get
\[
s-t+1 = a\left ( \frac{I:g}{I} (-1) \right )
\]
where $a(M) := \min \{ j \in {\mathbb Z} \ | \ [M]_j \neq 0 \}$ denotes the
initial degree.  Since $s \geq 2t$ we conclude
\[
a\left (\frac{I:g}{I} \right ) = s-t \geq \frac{s}{2} = \frac{\reg (R/I)}{2} =
\frac{\reg I -1}{2}.
\]

Conversely, if $I$ has the subspace property then the exact sequence
(\ref{exseq}) gives injectivity for the multiplication by $g$ on $R/I$ in
degrees $\leq s/2$, and the self-duality of $R/I$ (up to shift) gives the
surjectivity for the second half.
\end{proof}

\begin{lemma} \label{sub-pres}
Suppose that $I$ has the subspace property with respect to $g$.  Let $\ell$
be a general linear form.  If $\depth R/I > 0$ then $I+\ell R$ has the
subspace property with respect to $g$.
\end{lemma}

\begin{proof}
Consider the commutative diagram
\[
\begin{array}{ccccccccccccccccccc}
&& 0 && 0 && 0 && 0 \\
&& \downarrow && \downarrow && \downarrow && \downarrow \\
0 & \rightarrow & \frac{I:g}{I} (-2) & \rightarrow & \frac{R}{I} (-2) &
\stackrel{g}{\rightarrow} & \frac{R}{I} (-1) & \rightarrow & \frac{R}{I+gR}
(-1) & \rightarrow & 0 \\
&& \phantom{\ell} \downarrow \ell &&  \phantom{\ell} \downarrow \ell &&
 \phantom{\ell} \downarrow \ell &&  \phantom{\ell} \downarrow \ell \\
0 & \rightarrow & \frac{I:g}{I} (-1) & \rightarrow & \frac{R}{I} (-1) &
\stackrel{g}{\rightarrow} & \frac{R}{I}  & \rightarrow & \frac{R}{I+gR}
 & \rightarrow & 0 \\
&& \downarrow && \downarrow && \downarrow && \downarrow \\
0 & \rightarrow & \frac{(I+\ell R):g}{I+\ell R} (-1) & \rightarrow &
\frac{R}{I+\ell R} (-1) & \stackrel{g}{\rightarrow} & \frac{R}{I+\ell R} &
\rightarrow & \frac{R}{I+(g,\ell)R} & \rightarrow & 0
\end{array}
\]
It follows that
\[
a\left ( \frac{(I+\ell R):g}{I+\ell R} \right ) = a \left ( \frac{I:g}{I}
\right ) \geq \frac{\reg I -1}{2} = \frac{\reg (I+\ell R) -1}{2}.
\]
\end{proof}

In particular, we see that if $X$ is an \acm subscheme with the subspace
property then its general hyperplane section also has this property, provided
that $\dim X \geq 1$.  In view of the last part of Definition \ref{WSP}, where
we extended the weak Lefschetz property to the non-Artinian case, we have the
following corollary.

\begin{corollary}
Suppose that an ideal $I$ has the subspace property.  Then $I$ has the weak
Lefschetz property.
\end{corollary}

\begin{proof}
It follows from Lemma \ref{artin-sub-iff-wlp} and Lemma \ref{sub-pres}.
\end{proof}

\begin{lemma}
The ideals constructed in the proof of Theorem \ref{maingorthm} have the
subspace property.
\end{lemma}

\begin{proof}

In the proof of Theorem \ref{maingorthm} we use $Z \subset Z_{c-1,t} \subset
G_{c-1,s+1,t}$, where $s \geq 2t$.  Thus we get for $g := L_{s+1-t+ \lfloor
\frac{c-3}{2} \rfloor} = L_{s-t+\lfloor \frac{c-1}{2} \rfloor}$ that
\[
I_{G_{c-1,s+1,t}} : g = I_{G_{c-1,s,t}} =: K' \hbox{\hskip 1cm and \hskip
1cm} I_Z : g = I_Z
\]
(this follows from Theorem \ref{max-gor} and Theorem \ref{acm-scheme}).

Put $K := I_{G_{c-1,s+1,t}}$, $I := I_Z$, $J := K:I$.  Consider the
commutative diagram
\[
\begin{array}{cccccccccccccccc}
&&0 &&0 &&0 \\
&&\downarrow &&\downarrow &&\downarrow \\
0 & \rightarrow & \frac{R}{K:g}(-1) & \rightarrow & \left ( \frac{R}{I}
\oplus \frac{R}{J:g} \right ) (-1) &
\rightarrow & \left ( \frac{R}{(I+J):g} \right ) (-1) & \rightarrow & 0 \\
&&\downarrow &&\downarrow &&\downarrow \\
0 & \rightarrow & \frac{R}{K} & \rightarrow & \frac{R}{I} \oplus \frac{R}{J}
& \rightarrow & \frac{R}{I+J} & \rightarrow & 0 \\
&&\downarrow &&\downarrow &&\downarrow \\
0 & \rightarrow & \frac{R}{K+gR} & \rightarrow & \frac{R}{I+gR} \oplus
\frac{R}{J+gR} & \rightarrow & \frac{R}{I+J+gR} & \rightarrow & 0 \\
&&\downarrow &&\downarrow &&\downarrow \\
&& 0 && 0 && 0
\end{array}
\]
Since $I$ and $J:g$ are linked by $K' = K:g$, the ideal $(I+J):g = I+(J:g)$
is Gorenstein of codimension $c$ and has the $h$-vector
$(1,c,h_2,\dots,h_t,\dots,h_t, h_{t-1},\dots,c,1)$ with the final 1 occurring
in degree $s-1$.  Hence the rightmost column shows that $I+J+gR$ is a perfect
ideal of codimension $c$ in $R$ and
\[
a\left ( \frac{(I+J):g}{I+J} \right ) = s-t \geq \frac{s}{2} = \frac{\reg
(I+J) -1}{2}
\]
Therefore $I+J$ has the subspace property with respect to $g$.
\end{proof}

\begin{remark}\label{compare-harima-2}
The fact that our arithmetically Gorenstein schemes possess the subspace
property shows that in some sense the approach of Harima \cite{harima}
has an analog in higher dimension.  We find a large subset of $G_{c,s,t}$
lying on a hyperplane, and after it is removed we are left with
$G_{c,s-1,t}$, which is the higher dimensional analog of his ideal quotient
trick.   See Remark \ref{compare-harima-1}.

 On the other hand, the choice of the hyperplane which gives the subspace
property for a subscheme $X$ can be counter-intuitive, and in particular it is
not necessarily the hyperplane containing the ``largest'' subset of $X$, as is
shown in the following example.

Let $X$ consist of six points in $\proj{2}$, three of which are generically
chosen on a line $\lambda_1$, and three more generally chosen points on a
line $\lambda_2$.  Note that $X$ is a complete intersection, and $\reg(I_X) =
4$.  Let $\ell$ be a linear form.  We claim that if $\ell$ is the form
defining $\lambda_1$ or $\lambda_2$ then it does not give the subspace
property, but that if $\ell$ vanishes on one point of $X$ that lies on
$\lambda_1$ and one point of $X$ that lies on $\lambda_2$ then $\ell$ does
give the subspace property.  Indeed, in all three cases $(I+\ell R)/\ell R$ is
perfect in $R/\ell R$, but in the first two cases the initial degree of $(I:_R
\ell) /I$ is 1 (violating the second condition of the definition), while in
the last case this initial degree is 2, satisfying the condition.
\end{remark}


\section{Extremal Graded Betti Numbers}

The goal of this section is to show that the \aG schemes constructed in
Section \ref{tie-together} have the maximal graded Betti numbers (among the
class of \aG schemes with the same Hilbert function and the Weak Lefschetz
Property).  To this end we first compute the Betti numbers of the sum of two
geometrically (Gorenstein) linked ideals.  Second, in order to get upper
bounds for the Betti numbers we compare the Betti numbers of a scheme with
those of a hyperplane section.  This is done in the more general framework
of modules, i.e.\ we compare Betti numbers of $M$ and $M/\ell M$ where $M$ is
a graded module and $\ell$ is any linear form.

Let $X,Y \subset \proj{n}$ be geometrically linked subschemes of codimension
$c \leq n$.  We have already seen that the Hilbert function of $X \cap Y$ is
determined by the Hilbert function of $X$ alone.  The same is true of the
graded Betti numbers under an additional hypothesis.

\begin{lemma}
If $X,Y \subset \proj{n}$ are geometrically linked and $\reg (X \cup Y) \geq 2
\cdot \reg (X)$ then (for all $i$)
\[
\TT_i^R (R/(I_X+I_Y),K) \cong \TT_i^R(R/I_X,K) \oplus \TT_{i-1}^R
(K_X,K)(-\reg(X\cup Y)-c+n+2)
\]
where $K_X := \Ext_R^c (R/I_X,R)(-n-1)$ denotes the canonical module of $X$.
\end{lemma}

\begin{proof}
Following \cite{SV}, we define for a Noetherian graded module $M$,
\[
\begin{array}{rcl}
r(M) & := & \inf \{ n \in {\mathbb Z}\ | \ p_M(m) = h_M(m) \hbox{ for all } m
\geq n \} \\
e(M) & := & \sup \{ n \in {\mathbb Z}\ | \ [M]_n \neq 0 \},
\end{array}
\]
where $p_M(t)$ is the Hilbert polynomial of $M$.  In \cite{SV}, $r(M)$ is
called the {\em index of regularity} of
$M$. Since $X \cup Y$ is arithmetically Gorenstein by assumption, we get for
its index of regularity
\[
\begin{array}{rcl}
r(R/I_{X \cup Y}) & = & e(H^{n+1-c}_{\mathfrak m} (R/I_{X \cup Y})) +1 \\
& = & \reg (R/I_{X \cup Y}) - (n+1-c) +1 \\
& = & \reg (X \cup Y) + c-n-1
\end{array}
\]
Thus we obtain (cf.\ for instance \cite{nagel}, Lemma 2.5)
\[
K_X (-\reg(X \cup Y)-c+n+2) \cong I_Y/(I_X \cap I_Y) \cong (I_X+I_Y)/I_X
\]
and the exact sequence
\[
\begin{array}{cccccccccc}
&& 0 && 0 \\
&& \downarrow &&\downarrow \\
&& G_s && F_t \\
&& \downarrow &&\downarrow \\
&& \vdots && \vdots \\
&& \downarrow &&\downarrow \\
&& G_1 && F_1 \\
&& \downarrow &&\downarrow \\
&& G_0 && R \\
&& \downarrow &&\downarrow \\
0 & \rightarrow & K_X(-\reg(X\cup Y)-c+n+2) & \rightarrow & R/I_X & \rightarrow
& R/(I_X+I_Y) & \rightarrow & 0 \\
&& \downarrow &&\downarrow \\
&& 0 && 0
\end{array}
\]
(note that $X$ is not assumed to be arithmetically Cohen-Macaulay), where the
vertical sequences denote the corresponding minimal free resolutions.

Now we observe $K_X \cong H_{\goth m}^{n+1-c} (R/I_X)^\vee$ and
\[ e(H_{\goth m}^{n+1-c} (R/I_X)) +n+1-c \leq \reg(R/I_X) = \reg(X)-1.
\]
It follows that
\[
a(K_X) = -e(H_{\goth m}^{n+1-c} (R/I_X)) \geq -\reg(X)-c+n+2,
\]
Thus
\[
a(K_X(-\reg (X\cup Y) -c+n+2)) \geq \reg(X \cup Y) -\reg (X) \geq \reg (X)
\]
where the last inequality is by assumption.  We conclude that
$a(G_i) \geq \reg (X)+i$ for all $i \geq 0$.

On the other hand, it is well-known that the modules $F_i$ are generated in
degrees $\leq \reg (R/I_X) +i = \reg (X)-1+i$.  It follows that the mapping
cone procedure provides a free resolution of $R/I_X+I_Y$ which is minimal
because cancellation cannot occur.  This proves the claim.
\end{proof}

If $X$ is arithmetically Cohen-Macaulay, the graded Betti numbers of $K_X$ can
be computed from the graded Betti numbers of $X$.  Thus as a corollary we
obtain Theorem 7.1 of \cite{GHS1}:

\begin{corollary} \label{sum-betti}
If $X,Y \subset \proj{n}$ are geometrically linked \acm schemes such that
$\reg (X \cup Y) \geq 2 \cdot \reg (X)$ then
\[
\TT_i^R (R/(I_X+I_Y),K) \cong \TT_i^R (R/I_X ,K) \oplus
\TT_{c-i+1}^R (R/I_X,K)^\vee (-\reg(X\cup Y) -c+1).
\]
\end{corollary}

\begin{proof}
Since $R/I_X$ is Cohen-Macaulay, we get by applying $\Hom_R (-,R)$ to a
minimal free resolution of $R/I_X$ a minimal free resolution of $\Ext_R^c
(R/I_X ,R)$, i.e.
\[
\TT_i^R (\Ext_R^c (R/I_X,R),K) \cong \TT_{c-i}^R (R/I_X ,K)^\vee.
\]
Thus the Lemma implies our assertion.
\end{proof}

\begin{lemma} \label{tor-seq}
Let $M$ be a graded $R$-module, $\ell \in R$ a linear form.  Then there is an
exact sequence of graded $R$-modules (where $\bar R := R/ \ell R)$:
\[
\begin{array}{c}
\cdots \rightarrow \TT_{i-1}^{\bar R} ((0:_M \ell),K)(-1) \rightarrow
\TT_i^R (M,K) \rightarrow \TT_i^{\bar R} (M/ \ell M ,K) \rightarrow \cdots
\hbox{\hskip 1in} \\
\hfill \cdots \rightarrow \TT_1^R (M,K) \rightarrow \TT_1^{\bar R} (M/ \ell
M,K)  \rightarrow 0.
\end{array}
\]
\end{lemma}

\begin{proof}
Consider the exact sequence
\[
0 \rightarrow R(-1) \stackrel{\ell}{\rightarrow} R \rightarrow \bar R
\rightarrow 0.
\]
Tensoring with $M$ and taking homology provides $\TT_i^R (M,\bar R) = 0$ for
all $i \geq 2$ and
\[
\TT_1^R (M,\bar R) \cong \ker \left ( M(-1) \stackrel{\ell}{\rightarrow} M
\right ) = (0:_M \ell)(-1).
\]
Now we compute $\TT_i^R (M,\bar R)$ using a minimal free resolution of $M$:

\medskip
\noindent
$ \displaystyle
{\bf F}_\bullet \hskip 1.5in 0 \rightarrow F_s \rightarrow F_{s-1} \rightarrow
\cdots \rightarrow F_0 \rightarrow M \rightarrow 0.
$

\medskip

\noindent Tensoring by $\bar R$ gives the complex (with $\bar F_i := F_i
\otimes_R \bar R$)

\medskip
\noindent
$ \displaystyle
{\bf F}_\bullet \otimes_R \bar R
\hskip .8in 0 \rightarrow \bar F_s \rightarrow \cdots
\rightarrow
\bar F_2 \stackrel{\alpha}{\longrightarrow} \bar F_1
\stackrel{\beta}{\longrightarrow}
\bar F_0 \rightarrow M/ \ell M \rightarrow 0.
$

\medskip

\noindent Its homology is
\[
H_i ({\bf F}_\bullet \otimes \bar R ) =
\left \{
\begin{array}{ll}
0 & \hbox{if $i \leq 0$} \\
\TT_i^R (M,\bar R) & \hbox{if $i \geq 1$}
\end{array}
\right.
\cong
\left \{
\begin{array}{ll}
0 & \hbox{if $i \neq 1$} \\
(0:_M \ell)(-1) & \hbox{if $i = 1$}.
\end{array}
\right.
\]
Since ${\bf F}_\bullet$ is a minimal free resolution, the maps in ${\bf
F}_\bullet \otimes \bar R$ are minimal maps too.  Thus we obtain a minimal
free resolution of $\bar F_1 /\im \alpha$ as an $\bar R$-module:
\[
0 \rightarrow \bar F_s \rightarrow \cdots
\rightarrow
\bar F_2 \stackrel{\alpha}{\longrightarrow} \bar F_1 \rightarrow \bar F_1 /\im
\alpha \rightarrow 0.
\]
It follows that
\begin{equation}\label{number1}
\TT_i^{\bar R} (\bar F_1 / \im \alpha ,K) \cong \bar F_{i+1} \otimes_{\bar R} K
\cong \TT_{i+1}^R (M,K) \hskip .5in \hbox{for all $i \geq 0$}.
\end{equation}
The exact sequence
\[
0 \rightarrow \bar F_1 /\ker \beta \rightarrow \bar F_0 \rightarrow M/ \ell M
\rightarrow 0
\]
implies
\begin{equation}\label{number2}
\TT_i^{\bar R} (\bar F_1 / \ker \beta,K) \cong \TT_{i+1}^{\bar R} (M/\ell M,K)
\hskip .5 in \hbox{for all $i \geq 0$}.
\end{equation}
Using $\ker \beta / \im \alpha = H_1 ({\bf F}_\bullet \otimes \bar R) \cong
(0:_M \ell )(-1)$ we obtain the exact sequence
\[
0 \rightarrow (0:_M \ell)(-1) \rightarrow \bar F_1 /\im \alpha \rightarrow
\bar F_1 /\ker \beta \rightarrow 0.
\]
The associated long exact Tor sequence reads as
\[
\begin{array}{c}
\cdots \rightarrow \TT_1^{\bar R} (\bar F_1 /\ker \beta, K) \rightarrow
\TT_0^{\bar R} ((0:_M \ell),K)(-1) \rightarrow \TT_0^{\bar R} (\bar F_1 /\im
\alpha,K) \hskip 1in \\
\hfill \rightarrow \TT_0^{\bar R} (\bar F_1/\ker \beta,K)
\rightarrow 0.
\end{array}
\]
Taking into account the isomorphisms (\ref{number1}) and (\ref{number2}), our
claim follows.
\end{proof}

\begin{notation}
We denote $\displaystyle \left [ \tor_i^R (M,K) \right ]_j := \rank_K \left [
\TT_i^R (M,K) \right ]_j$.
\end{notation}

\begin{corollary} \label{lastcor}
With the notation of Lemma \ref{tor-seq} we have for all $i \geq 1$ and all $j
\in {\mathbb Z}$:
\[
\left [ \tor_i^R (M,K) \right ]_j \leq \left [ \tor_i^{\bar R} (M/\ell M,K)
\right ]_j + \left [ \tor_{i-1}^{\bar R} ((0:_M \ell),K) \right ]_{j-1}.
\]
Furthermore,
\[
\begin{array}{rl}
\left [ \TT_i^R (M,K) \right ]_j \cong \left [ \TT_i^{\bar R} (M/\ell M,K)
\right ]_j & \hbox{if $i=1$ and $j \leq a+i-1$} \\
& \hbox{or $i \geq 2$ and $j \leq a+i-2$}
\end{array}
\]
where $a := a(0:_M \ell)$.
\end{corollary}

\begin{proof}
Using $[\TT_i^{\bar R} ((0:_M \ell),K)]_j = 0$ if $j < a+i$ the claim follows
by analyzing the sequence given in Lemma \ref{tor-seq}.
\end{proof}

\begin{remark}
In the special case where $\ell$ is a non-zero divisor for $M$ we have
$a(0:_M \ell)$ $ = \infty$.  Then we get back the well-known fact that graded
Betti numbers do not change under such hyperplane sections.
\end{remark}

Now we specialize the previous results to Gorenstein $K$-algebras.

\begin{proposition}\label{toreq}
Let $A = R/I$ be a graded Artinian Gorenstein $K$-algebra and let $c =
 \rank_K[R]_1$.  Let $\ell \in R$ be any linear form.  Put $s :=
\reg (A)$ and $a :=$ $a(0:_A \ell)$.  Then we have for all $i \in {\mathbb Z}$
\[
\begin{array}{l}
\left [ \tor_i^R (A,K) \right ]_j = \\ \\

\hskip 1cm
\left \{
\begin{array}{ll}
\left [ \tor_i^{\bar R} (A/\ell A,K) \right ]_j
& \hbox{if } j \leq a+i-2 \\
\leq \left [ \tor_i^{\bar R} (A/\ell A,K) \right ]_j +
\left [ \tor_{c-i}^{\bar R} (A/\ell A,K) \right ]_{s+c-j}
& \hbox{if } a+i-1 \leq j \leq s-a+i+1 \\
\left [ \tor_{c-i}^{\bar R} (A/\ell A,K) \right ]_{s+c-j}
& \hbox{if } j \geq s-a+i+2
\end{array}
\right.
\end{array}
\]
\end{proposition}

\begin{proof}
We denote by $^\vee := \Hom_R (-,K)$ the dualizing functor with respect to $K
\cong R/{\goth m}$.  It is exact.  Thus, the exact sequence
\[
0 \rightarrow (0:_A \ell)(-1) \rightarrow A(-1) \rightarrow A \rightarrow
A/\ell A \rightarrow 0
\]
provides the exact sequence
\[
0 \rightarrow (A/\ell A)^\vee (-1) \rightarrow A^\vee (-1) \rightarrow A^\vee
\rightarrow (0:_A \ell)^\vee \rightarrow 0.
\]
Since $A$ is Artinian and Gorenstein, we have $A^\vee \cong A(s)$.  Consider
the commutative diagram
\[
\begin{array}{ccccccccccccc}
0 & \rightarrow & (0:_A \ell)(-1) & \rightarrow & A(-1) &
\stackrel{\ell}{\longrightarrow} & A & \rightarrow & A/\ell A & \rightarrow & 0
\\ &&&& \phantom{\wr} \downarrow \wr && \phantom{\wr} \downarrow \wr \\
0 & \rightarrow & (A/\ell A)^\vee (-s-1) & \rightarrow & A^\vee (-s-1) &
\stackrel{\ell}{\longrightarrow} & A^\vee (-s) & \rightarrow & (0:_A \ell)^\vee
(-s)  & \rightarrow & 0.
\end{array}
\]
We conclude that
\begin{equation} \label{star}
0:_A \ell \cong (A/\ell A)^\vee (-s).
\end{equation}
Put $\bar R = R / \ell R$.  Since $A/\ell A$ is Artinian we get by local
duality
\[
(A/\ell A)^\vee \cong H^0_{\goth m} (A/\ell A)^\vee \cong \Ext_{\bar
R}^{c-1} (A/\ell A,\bar R)(-c).
\]
It follows that
\[
\TT_i^{\bar R} ((A/\ell A)^\vee,K) \cong \left ( \TT_{c-1-i}^{\bar R} (A/\ell
A,K) \right )^\vee (-c)
\]
and therefore
\[
\TT_i^{\bar R} (0:_A \ell,K) \cong \left ( \TT_{c-1-i}^{\bar R} (A/\ell A,K)
\right )^\vee (-s-c) \hskip 1cm \hbox{(as $\bar R$-modules)}.
\]
In particular, we get for all $i,j \in {\mathbb Z}$ isomorphisms of $K$-vector
spaces
\[
\left [
\TT_i^{\bar R} (0:_A \ell,K) \right ]_j \cong
\left [ \TT_{c-1-i}^{\bar R} (A/\ell A,K) \right ]_{s+c-j}.
\]
Hence Corollary \ref{lastcor} proves the claim for all pairs $(i,j)$ with $j
\leq s-a+i+1$.  Since the minimal free resolution of $A$ is self-dual we have
\[
\TT_i^R (A,K) \cong \left ( \TT_{c-i}^R (A,K) \right )^\vee (-s-c).
\]
Using what we have already proved we get for $j \geq s-a+i+2$
\[
\left [ \tor_i^R (A,K) \right ]_j = \left [ \tor_{c-i}^R (A,K) \right ]_{s+c-j}
= \left [ \tor_{c-i}^{\bar R} (A/\ell A,K) \right ]_{s+c-j}
\]
concluding the proof.
\end{proof}

\begin{remark} Proposition \ref{toreq} provides a somewhat quantitative
version of the following fact: The more general $\ell$ is, the more the graded
Betti numbers of $A$ as an $R$-module are determined by the graded Betti
numbers of $A/\ell A$ as an $R/\ell R$-module.  
\end{remark}

If $\reg(A) \leq 2 \cdot
a(0:_A \ell) -3$, the formula in Proposition \ref{toreq} simplifies as
follows:

\begin{corollary} \label{fax-8.9}
Using the notation and assumptions of Proposition \ref{toreq}, suppose in
addition that $s \leq 2a-3$.  Then we have for all $i \in {\mathbb Z}$
\[
[ \tor_i^R (A,K)]_j = 
\left \{
\begin{array}{ll}
[\tor_i^{\bar R} (A/ \ell A,K)]_j & \hbox{if } j \leq a+i-2 \\  \\
\left [ \tor_{c-i}^{\bar R} (A/ \ell A,K) \right ]_{s+c-j} & \hbox{if } j \geq
s-a+i+2.
\end{array}
\right.
\]
\end{corollary}

\begin{proof}
If $s \leq 2a-3$ then there is no integer $j$ satisfying $a+i-1 \leq j \leq
s-a+i+1$.  Thus Proposition \ref{toreq} proves the claim.
\end{proof}

\begin{notation}
Let $\underline{h} = (1,h_1,\dots,h_s)$ be the $h$-vector of an Artinian
$K$-algebra. Let $c \geq h_1$ be an integer.  Then there is a uniquely
determined lex-segment ideal $I \subset K[z_1,\dots,z_{c}] =: T$ such that
$T/I$ has
$\underline{h}$ as its Hilbert function.  We define
\[
\beta_{i,j} (\underline{h},c) := \left [ \tor_i^T (T/I,K) \right ]_j.
\]
If $c=h_1$ we simply write $\beta_{i,j}(\underline{h})$ instead of
$\beta_{i,j}(\underline{h},h_1)$.
\end{notation}

\begin{remark}
The numbers $\beta_{i,j}(\underline{h},c)$ can be computed numerically
without considering lex-segment ideals.  Explicit formulas can be found in 
\cite{EK}.
\end{remark}

\begin{theorem}[\cite{bigatti}, \cite{hulett}] \label{big-hul}
If $I \subset R = K[x_0,\dots,x_n]$ is a Cohen-Macaulay ideal of codimension
$c$ and $\underline{h}$ is the $h$-vector of $R/I$ then we have for all $i,j
\in {\mathbb Z}$
\[
\left [ \tor_i^R (R/I,K) \right ]_j \leq \beta_{i,j} (\underline{h},c).
\]
\end{theorem}

We are now ready for the main result of this paper.

\begin{theorem} \label{main-thm-paper}
Let $\underline{h} = (1,h_1,h_2,\dots,h_t,\dots, h_s)$ be an SI-sequence where
$h_{t-1 } < h_t = \cdots = h_{s-t} > h_{s-t+1}$.  Put $\underline{g} =
(1,h_1-1,h_2-h_1,\dots,h_t-h_{t-1})$.  Then we have
\begin{itemize}
\item[(a)] If $A = R/I$ is a Gorenstein $K$-algebra with $c = \codim I$ and
having  $\underline{h}$ as $h$-vector and an Artinian reduction which has the
weak Lefschetz property, then
\[
\left [ \tor_i^R (A,K) \right ]_j \leq
\left \{
\begin{array}{ll}
\beta_{i,j} (\underline{g},c-1) & \hbox{if $j \leq s-t+i-1$} \\
\beta_{i,j} (\underline{g},c-1) + \beta_{c-i,s+c-j} (\underline{g},c-1) &
\hbox{if $s-t+i \leq j \leq t+i$} \\
\beta_{c-i,s+c-j} (\underline{g},c-1) & \hbox{if $j \geq t+i+1$}
\end{array}
\right.
\]
\item[(b)] Suppose that $K$ has sufficiently many elements. 
Then there is a reduced, non-degenerate
arithmetically Gorenstein subscheme $X \subset \proj{n} = \hbox{Proj} (R)$ of
codimension $c$ with the subspace property (and hence the Artinian reduction
of $A = R/I_X$ has the weak Lefschetz property) and $h$-vector
$\underline{h}$, and with equality holding for all $i,j \in {\mathbb Z}$ in
{\rm (a)}.  Indeed, we can take $X = G_c(\underline{h})$ whenever the set
${\mathcal N}_{c-1,s+1,t}$ is chosen sufficiently general, i.e.\ satisfies
the condition given in Theorem \ref{maingorthm}.
\end{itemize}
\end{theorem}

\begin{proof}
We first prove (a).  We may assume that $A$ is Artinian.  Then it has the weak
Lefschetz property, i.e.\ there is an $\ell \in [R]_1$ such that the
multiplication map
\[
[A]_{i-1} \stackrel{\ell}{\longrightarrow}
[A]_i
\]
 has maximal rank for all $i \in {\mathbb Z}$.  It follows that the Hilbert
function of $A/\ell A$ is $\underline{g}$.  Thus Theorem \ref{big-hul} gives
\[
\left [ \tor_i^{\bar R} (A/\ell A,K) \right ]_j \leq \beta_{i,j}
(\underline{g},c-1) 
\]
where $\bar{R} = R/l R$. 
Therefore Proposition \ref{toreq} provides for all $i,j \in {\mathbb Z}$
\[
\left [ \tor_i^R (A,K) \right ]_j \leq \beta_{i,j}(\underline{g},c-1) +
\beta_{c-i,s+c-j} (\underline{g},c-1).
\]
The claim then follows because $\beta_{i,j} (\underline{g},c-1) = 0$ if $j >
t+i$.

For (b), according to Theorem \ref{maingorthm} and Lemma
\ref{sub-pres} it suffices to compute the graded Betti numbers of the
subscheme $X = G_c(\underline{h}) \subset
\proj{n}$ constructed in  Theorem \ref{maingorthm}.  The starting
point was an \acm subscheme $Z \subset \proj{n} = \hbox{Proj} (R)$ of
codimension $c-1$  with
$h$-vector $\underline{g}$ which satisfies
\[
\left [ \tor_i^R (R/I_Z ,K) \right ]_j =
\beta_{i,j}(\underline{g},c-1) \ \ \hbox{for all $i,j \in {\mathbb Z}$}
\]
thanks to Corollary \ref{max-betti}.  Let $Y$ be the
residual to $Z$ in $G_{c-1,s+1,t}$.  Then we have $I_X = I_Z + I_Y$.  The
graded Betti numbers were computed in Corollary \ref{sum-betti}:
\[
\begin{array}{rcl}
\left [ \tor_i^R (R/I_X,K) \right ]_j & = & \left [ \tor_i^R (R/I_Z,K)
\right ]_j + \left [ \tor_{c-i}^R (R/I_Z,K) \right ]_{s+c-j} 
 \\ \\
& = & \beta_{i,j} (\underline{g},c-1) + \beta_{c-i,s+c-j} (\underline{g},c-1).
\end{array}
\]
Using again the fact that $\beta_{i,j} (\underline{g},c-1) = 0$ if $j > t+i$,
we get the result.
\end{proof}

Next we prove that if $I \subset R$ is a Gorenstein ideal of codimension $c$
such that $R/I$ has the Weak Lefschetz property and maximal $h$-vector
$\underline{h}$ with respect to $c,t,s$, then the graded Betti numbers of
$R/I$ are determined by $c,t,s$ (and equal to the ones of
$R/J_c(\underline{h})$).

\begin{corollary} \label{fax-8.14}
Let $c,s,t$ be positive integers, where either $s = 2t$ or $s \geq 2t+2$. 
Define $\underline{h} = (h_0,\dots,h_s)$ by 
\[
h_i = \left \{
\begin{array}{rl}
\binom{c-1+i}{c-1} & \hbox{if $0 \leq i \leq t$}; \\
\binom{c-1+t}{c-1} & \hbox{if $t \leq i \leq s-t$};\\
\binom{s-i+c-1}{c-1} & \hbox{if $s-t \leq i \leq s$}
\end{array}
\right.
\]
Let $I \subset R$ be a Gorenstein ideal of codimension $c$ such that $R/I$
has the weak Lefschetz property and $h$-vector $\underline{h}$.  Then the
minimal free resolution of $R/I$ has the shape
\[
0 \rightarrow R(-s-c) \rightarrow F_{c-1} \rightarrow \cdots \rightarrow F_1 
\rightarrow R \rightarrow R/I \rightarrow 0
\]
where 
\[
F_i = R(-t-i)^{\alpha_i} \oplus R(t-s-i)^{\gamma_i} \ \ \hbox{ and } \ \
\alpha_i = \binom{c+t-1}{i+t} \binom{t-1+i}{t} = \gamma_{c-i}.
\]
\end{corollary}

\begin{proof}
We may assume that $R/I$ is Artinian and has the weak Lefschetz property with
respect to $\ell = x_{c-1}$.  Thus $A/\ell A \cong \bar
R/(x_0,\dots,x_{c-2})^{t+1}$.  Its graded Betti numbers are easily
computed. For the reader's convenience we include the argument. 

Let $\ffi: F = \bar{R}^{c+t-1} \to G = \bar{R}^{t+1}$ be the graded
homomorphism given by the matrix 
$$
M = \left ( \begin{array}{ccccccc}
x_0 & x_1 & \ldots & x_{c-2} & & &0 \\
& x_0 & x_1 & \ldots & x_{c-2} & \\
& & \ddots & & & \ddots \\
0 & & & x_0 & x_1 & \ldots & x_{c-2} 
\end{array} \right ) \in \bar{R}^{t+1, c+t}. 
$$
The ideal of maximal minors of $M$ is $(x_0,\dots,x_{c-2})^{t+1}$ and it is
resolved by the  Eagon-Northcott complex 
$$
\begin{array}{c}
0 \to \wedge^{c+t-1} F \otimes S_{c-2}(G)  \to \wedge^{c+t-2}
F \otimes S_{c-3}(G)  \to \ldots \hbox{\hskip 1in} \\
\hfill \to \wedge^{t+1} F \otimes S_0(G) \to (x_0,\dots,x_{c-2})^{t+1} \to 0
\end{array}
$$
In particular,
we get for all $i,j \in {\mathbb Z}$:
\[
\left [ \tor_i^{\bar R} (A/\ell A,K) \right ]_j = \beta_{i,j}(\underline{g}),
\ \ \hbox{ where } \underline{g} = (g_0,\hdots,g_t) \ \ \hbox{ and } g_i =
\binom{c-2+i}{i}. 
\]
If $s \geq 2t+3$ the claim follows by Corollary \ref{fax-8.9} because  
$a := a(0:_A \ell) = s-t$ which then forces $s \leq 2a-3$.  If $s=2t+2$ then 
we  apply Proposition \ref{toreq}. It proves our assertion since we have for
$i > 0$ that 
$ [ \tor_i^{\bar R} (A/\ell A,K)  ]_j = 0$ if $j \neq i+t$.  
If $s=2t$ then Theorem \ref{main-thm-paper} shows that
$[\tor_i^R(A,K)]_j = 0$ unless $(i,j) = (0,0)$, $(i,j) =
(c,s+c)$ or $j = i+t$ with $1 \leq i \leq c-1$.  But then
$[\tor_i^R(A,K)]_{i+t}$ can be computed recursively from the 
Hilbert function of $A$.  (A similar computation can be found on page 4386 of
\cite{uwetrans}.)
\end{proof}

\begin{example}
Corollary \ref{fax-8.14} cannot be extended to the case $s=2t+1$.  Indeed,
consider the case $c=4$, $s=5$ and $t=2$.  Then 
\[
\underline{h} = (1,4,10,10,4,1).
\]
The Betti numbers predicted by Corollary \ref{fax-8.14} are represented by
the Macaulay \cite{macaulay} diagram
\begin{verbatim}
; total:      1    16    30    16     1 
; --------------------------------------
;     0:      1     -     -     -     - 
;     1:      -     -     -     -     - 
;     2:      -    10    15     6     - 
;     3:      -     6    15    10     - 
;     4:      -     -     -     -     - 
;     5:      -     -     -     -     1 
\end{verbatim}
However, Boij \cite{boijjpaa} has shown that the generic Betti numbers of a
Gorenstein Artin algebra with this Hilbert function are
\begin{verbatim}
; total:      1    10    18    10     1 
; --------------------------------------
;     0:      1     -     -     -     - 
;     1:      -     -     -     -     - 
;     2:      -    10     9     -     - 
;     3:      -     -     9    10     - 
;     4:      -     -     -     -     - 
;     5:      -     -     -     -     1 
\end{verbatim}
We do not know if there exists a reduced arithmetically Gorenstein scheme
with these Betti numbers, but it would be very surprising to us if it did not
exist.

Also, it was pointed out to us by T.\ Harima that in Example 4.4 of
\cite{ikeda}, H.\ Ikeda constructed an interesting example of an Artinian
Gorenstein algebra with the above $h$-vector, not having the weak Lefschetez
property.  It is not known whether this Artinian Gorenstein ideal can be
lifted to  a zero-dimensional reduced, arithmetically Gorenstein subscheme in
$\proj{4}$ or not. 
\end{example}


\section{Simplicial Polytopes}

If the dimension of  the ambient space is large enough we can choose the
linear forms in the sets ${\mathcal M}_{c,t}$ and ${\mathcal N}_{c,s,t}$ as
variables.  Then the resulting ideals in Theorems \ref{anyO-seq}  and
\ref{maingorthm} are monomial ideals.   Thus they 
correspond to simplicial complexes.  Relating our results to the
constructions of Billera and Lee in \cite{BL2}, we will show that the
simplicial complexes we are getting are very particular.  They are
triangulations of balls or polytopes and have extremal properties with
respect to their Betti numbers.

For the convenience of the reader we recall some basic definitions.  Further
details can be found in \cite{BH}, \cite{Hi} or \cite{St2}.

A simplicial complex on a vertex set $V = \{ v_0,\dots,v_n \}$ is a
collection $\Delta$ of subsets of $V$ such that all $\{v_i\}$ belong to
$\Delta$ and $F \in \Delta$ whenever $F \subset G$ for some $G \in \Delta$. 
The elements of $\Delta$ are called the {\em faces} of $\Delta$.  The faces
$\{ v_i \}$ are called {\em vertices}, and the maximal faces under inclusion
are called {\em facets}.  The dimension of a face $F \in \Delta$ is $\dim F
:= |F|-1$.  The dimension of $\Delta$ is defined to be $\dim \Delta = \max \{
\dim F | F \in \Delta \}$.  The complex is called {\em pure} if all its
facets have the same dimension.

Let $\Delta$ be a $(d-1)$-dimensional simplicial complex.  The number of
$i$-dimensional faces of $\Delta$ is denoted by $f_i$.  Putting $f_{-1} :=1$,
the vector of positive integers $(f_{-1},f_0,\dots,f_{d-1})$ is called
the {\em $f$-vector} of $\Delta$.  The complex $\Delta$ is called {\em
shellable} if $\Delta$ is pure and its maximal faces can be ordered
$F_1,\dots,F_p$ so that the complex $(\bigcup_{j=1}^{i-1} \bar F_j) \cap
\bar F_i$ is pure of dimension $d-2$ for every $i$ with $2 \leq i \leq p$. 
Here $\bar F_j$ is the complex
\begin{equation}\label{def-of-bar-F}
\bar F_j := \{ \sigma \in \Delta | \sigma \subset F_j \}.
\end{equation}
Let $\Delta$ be a simplicial complex on the vertex set $V =
\{v_0,\dots,v_n\}$.  The ideal $I_\Delta \subset R = K[x_0,\dots,x_n]$ which
is generated by all square-free monomials $x_{i_1}x_{i_2}\cdots x_{i_q}$ such
that $\{ v_{i_1},\dots,v_{i_q}\} \notin \Delta$, is called the {\em
Stanley-Reisner ideal of $\Delta$}.  The Stanley-Reisner ring $K[\Delta]$ of
$\Delta$ is the factor ring $R/I_\Delta$.  The Hilbert series
$H_{K[\Delta]}(z) := \sum_{j=0}^\infty h_{K[\Delta]}(j) \cdot z^j$ can be
written as
\[
H_{K[\Delta]}(z) = \frac{Q(z)}{(1-z)^d}
\]
where $d = \dim K[\Delta] = \dim \Delta+1$ and $Q(z) = \sum_{i=0}^t h_i z^i$
is a polynomial in ${\mathbb Q}[z]$.  The vector of integers
$(h_0,\dots,h_t)$ is called the {\em $h$-vector} of $\Delta$.  If $K[\Delta]$
is Cohen-Macaulay then the $h$-vector of $\Delta$ and the $h$-vector of
$K[\Delta]$ agree.  The $f$-vector and the $h$-vector of $\Delta$ are related
by the equality of polynomials
\[
\sum_{j=0}^t h_j z^j = \sum_{j=0}^d f_{j-1}\cdot z^j (1-z)^{d-j}.
\]
It follows in particular that $t \leq d = \dim \Delta +1$ and 
\[
n+1 = f_0 = h_1 + d.
\]
Moreover, one can compute the $h$-vector from the $f$-vector and vice versa. 
The Stanley-Reisner ideal $I_{\Delta}$ is determined by the facets of
$\Delta$.  In fact, one has (cf.\ for instance \cite{BH}, Theorem 5.1.4).

\begin{lemma} \label{fax-9.1}
The Stanley-Reisner ideal of a simplicial complex $\Delta$ is
\[
I_\Delta = \bigcap_F \mathfrak{p}_F
\]
where the intersection is taken over all facets of $\Delta$ and
$\mathfrak{p}_F$ denotes the prime ideal generated by all variables $x_i$
such that $v_i \notin F$.
\end{lemma}

The simplicial complex $\Delta$ is called Cohen-Macaulay if $K[\Delta]$ is a
Cohen-Macaulay ring.  Such a complex is pure.  A shellable simplicial
complex is Cohen-Macaulay (cf., for example, \cite{BH}, Theorem 5.1.13). 
Nevertheless, the possible $h$-vectors of shellable complexes are the same
as the $h$-vectors of Cohen-Macaulay algebras according to the following
result of Stanley \cite{St1} (cf.\ also \cite{BFS}).

\begin{theorem} \label{stanleysthm}
Let $\underline{h} = (h_0,\dots,h_t)$ be a sequence of positive integers. 
Then the following conditions are equivalent:
\begin{itemize}
\item[(a)] $\underline{h}$ is the $h$-vector of a shellable simplicial
complex.

\item[(b)] $\underline{h}$ is the $h$-vector of a Cohen-Macaulay algebra.

\item[(c)] $\underline{h}$ is an O-sequence.
\end{itemize}
\end{theorem}

We want to show that a similar phenomenon occurs for the maximal Betti
numbers of shellable simplicial
complexes. Let $\Delta$ be a simplicial complex on the vertex set $V = \{
v_0,\dots, v_n \}$.  The graded Betti numbers of $\Delta$ are defined to be
\[
\beta_{ij}(\Delta) := \left [ \tor_i^R (K[\Delta],K \right ]_j.
\]

\begin{remark} \label{rem-Betti-char} 
The Betti numbers of $K[\Delta]$ may depend on the characteristic of
$K$. In fact, Reisner exhibited in \cite{Reisner} the minimal triangulation
of the projective plane. Then $K[\Delta]$ is Cohen-Macaulay if and only if
$char (K) \neq 2$. This possible dependence on the characteristic is not
made explicit in the notation $\beta_{ij}(\Delta)$ but we will mention
assumptions on the characteristic if they are necessary. 
\end{remark} 

Let $\underline{h}$ be the $h$-vector of $\Delta$.  If $\Delta$ is
Cohen-Macaulay 
then it is known (cf.\ Theorem \ref{big-hul}) that we have for all $i,j \in
{\mathbb Z}$
\[
\beta_{ij}(\Delta) \leq \beta_{ij}(\underline{h}).
\]
The next result shows that this estimate is the best possible, even for
shellable complexes. 

\begin{proposition} \label{fax-9.3}
Let $\underline{h} = (h_0,\dots,h_t)$ be an O-sequence.  Then there is a
shellable simplicial complex $\Delta$ of dimension $2t-1$ with $h$-vector
$\underline{h}$ having maximal graded Betti numbers, i.e.\
$\beta_{ij}(\Delta) = \beta_{ij}(\underline{h})$ for $i,j \in {\mathbb Z}$.
\end{proposition}

\begin{proof}
Define the integer $n$ by $n+1 = c+2t$ where $c = h_1$.  We will use the
notation of Definition \ref{udef}.  We list the monomials in
LOIM$(\underline{h})$  in order, $m_1 <_r m_2 <_r \dots <_r m_p$.  Using the
bijection $V = \{ v_0,\dots,v_n \} \rightarrow \{u_1,\dots,u_{c+2t} \}$,
given by $v_i \mapsto u_{i+1}$, we define $F_i \subset V$ as the subset
corresponding to $\bar \beta_{c,t} (m_i)$.  Finally, we put
\[
\Delta := \bigcup_{i=1}^p \bar F_i
\]
where $\bar F_i$ is defined in (\ref{def-of-bar-F}).
This complex is shellable with shelling order $F_1,F_2,\dots,F_p$ according
to \cite{BL2}.

As the next step we define the set ${\mathcal M}_{c,t} \subset R =
K[x_0\,\dots,x_n]$ by
\[
{\mathcal M}_{c,t} := \left \{ M_0,\dots,M_{t+\lfloor \frac{c-1}{2}
\rfloor},L_0,\dots, L_{t+\lfloor \frac{c-2}{2} \rfloor} \right \},
\]
with $M_i := x_{2i}$ and $L_i := x_{2i+1}$.  We claim that
$I_{c,t}(\underline{h}) \subset R$ is the Stanley-Reisner ideal of
$\Delta$.  Indeed, the definitions immediately imply ${\mathfrak p}_{F_i} =
{\mathfrak p}_{c,t}(m_i)$.  But Lemma \ref{fax-9.1} provides $I_\Delta =
\bigcap_{i=1}^p \mathfrak{p}_{F_i}$ and Theorem \ref{decompofh} gives 
\[
I_{c,t}(\underline{h}) = \bigcap_{i=1}^p \mathfrak{p}_{c,t}(m_i).
\]
Thus we conclude that $I_\Delta = I_{c,t}(\underline{h})$.  Hence Theorem
\ref{anyO-seq} shows that $\underline{h}$ is the $h$-vector of $\Delta$, and
Corollary \ref{max-betti} gives $\beta_{ij}(\Delta) =
\beta_{ij}(\underline{h},c)$.  Finally, we get $\dim \Delta = n+1-h_1-1 =
2t-1$.
\end{proof}

Now we turn to polytopes.  A {\em $d$-polytope} $P$ is the $d$-dimensional
convex hull of a finite set of points in ${\mathbb R}^d$.  It is called {\em
simplicial} if all its proper faces are simplices.  Then the collection of
the proper faces of $P$ together with the empty set is called the {\em
boundary complex} $\Delta(P)$.  Identifying $j$-dimensional simplices with
their $j+1$ vertices, we consider $\Delta(P)$ as the $(d-1)$-dimensional
simplicial complex on the vertex set $V$, where $V = \{v_0,\dots,v_n \}$ is
the set of vertices of $P$.  The $h$-vector and the graded Betti numbers of
the polytope $P$ are defined to be the $h$-vector and the graded Betti
numbers of $\Delta(P)$, i.e. $\beta_{ij}(P) := \beta_{ij}(\Delta(P))$.

The famous $g$-theorem says that a sequence $\underline{h} = (h_0,\dots,h_s)$
of positive integers is the $h$-vector of a simplicial polytope if and only
if $\underline{h}$ is an SI-sequence.  This result was conjectured by
McMullen \cite{McM}.  Sufficiency of the condition was proved by Billera and
Lee in \cite{BL2} and the necessity by Stanley in \cite{stanley2}.

The main result of this section gives an optimal upper bound for the graded
Betti numbers of simplicial polytopes.

\begin{theorem} \label{fax-9.5} Suppose that $K$ is a field of
  characteristic zero. 
Let $\underline{h} = (h_0,\dots,h_s)$ be an SI-sequence where 
\[
h_{t-1} < h_t
= \cdots = h_{s-t} > h_{s-t+1}.
\]
 Put $\underline{g} = (1,h_1-h_0,\dots,h_t-h_{t-1})$.
Then we have
\begin{itemize}
\item[(a)] If $P$ is a simplicial $d$-polytope with $h$-vector
$\underline{h}$ then
\[
\beta_{ij}(P) \leq
\left \{
\begin{array}{ll}
\beta_{ij}(\underline{g}) & \hbox{if $j \leq s-t+i-1$}; \\
\beta_{ij}(\underline{g}) + \beta_{h_1-i,s+h_1-j}(\underline{g}) &
\hbox{if $s-t+i \leq j \leq t+i$}; \\
\beta_{h_1-i,s+h_1-j}(\underline{g}) & \hbox{if $j \geq t+i+1$}.
\end{array}
\right.
\]
(Observe that necessarily $d=s$.)

\item[(b)] The bounds given in {\rm (a)} are sharp: there is an
  $s$-polytope $P$ 
with $h$-vector $\underline{h}$ for which all of the bounds given in {\rm
  (a)}  are attained.
\end{itemize}
\end{theorem}

\begin{proof}
(a) It is well known that the boundary complex $\Delta(P)$ of $P$ is
Gorenstein (cf.\ for instance \cite{BH}, Corollary 5.5.6).  Moreover,
Stanley proves in \cite{stanley2} one direction of the $g$-theorem by
showing that $K[\Delta(P)]$ in fact has the Weak Lefschetz property.  Since
the Stanley-Reisner ideal $I_{\Delta(P)}$ has codimension $h_1$, we conclude
by Theorem \ref{main-thm-paper}.

We now turn to (b).  We will first describe the boundary complex following
Billera and Lee \cite{BL2} and then verify that it has the required
properties.  

Put $c=h_1$.  We are looking for a polytope with vertex set $V =
\{v_0,\dots,v_{s+c-1}\}$.  Let $V' := \{ v_0,\dots,v_{c-2+2t} \}$ and $V'' =
V \backslash V'$.  Now we proceed in a manner similar to our proof of
Proposition \ref{fax-9.3}.  We list the monomials in $LOIM(\underline{g})$
in order, $m_1 <_r m_2 <_r \cdots <_r m_p$.  Using the bijection $V'
\rightarrow \{ u_1,\dots,u_{c-1+2t} \}$, which sends $v_i \mapsto u_{i+1}$,
we define $F' \subset V$ as the subset corresponding to $\bar
\beta_{c-1,t}(m_i)$.  Put $F_i := F_i' \cup V''$ and define
\[
\Delta := \bigcup_{i=1}^p \bar F_i
\]
where again, $\bar F_i$ is defined in (\ref{def-of-bar-F}).  It follows
immediately that $\Delta$ is a pure $s$-dimensional simplicial complex.  Let
$\partial \Delta$ denote the pure $(s-1)$-dimensional simplicial complex on
the vertex set $V$ whose facets are the $(s-1)$-dimensional faces of $\Delta$
which are contained in exactly one facet of $\Delta$.  

Billera and Lee have shown in \cite{BL2}, section 7 (using the construction
of section 5) that $\partial \Delta$ is indeed the boundary complex of a
simplicial $s$-polytope.  In order to complete the proof we relate the
Stanley-Reisner ideals $I_\Delta$ and $I_{\partial \Delta}$ to the ideals
constructed in our sections 5 and 6, respectively.

Put $n = s+c-1$ and define the sets 
\[
{\mathcal N}_{c-1,s+1,t} = \left \{ M_0,\dots,M_{t+\lfloor
\frac{c-1}{2} \rfloor} ,L_0,\dots,L_{s-t+\lfloor \frac{c-1}{2} \rfloor}
\right \} \subset R = K[x_0,\dots,x_n]
\]
and
\[
{\mathcal M}_{c-1,t} = \left \{ M_0,\dots,M_{t+\lfloor \frac{c-2}{2} \rfloor}
,L_0,\dots,L_{t+\lfloor
\frac{c-3}{2} \rfloor} \right \}
\]
by $M_i := x_{2i}$ and 
\[
L_i := \left \{
\begin{array}{ll}
x_{2i+1} & \hbox{if $0 \leq i \leq t + \lfloor \frac{c-3}{2} \rfloor$} \\
x_{s+c-1-i} & \hbox{if $t+ \lfloor \frac{c-1}{2} \rfloor \leq i \leq
s-t+\lfloor \frac{c-1}{2} \rfloor$}.
\end{array}
\right.
\]
Thus, we have ${\mathcal N}_{c-1,s+1,t} = \{x_0,\dots,x_n \}$ and ${\mathcal
M}_{c-1,t} = \{ x_0,\dots,x_{c-2+2t} \}$.

Using the notation of Lemma \ref{fax-9.1}, we see that the ideals
$\mathfrak{p}_{F_i}$ and $\mathfrak{p}_{c,t}(m_i)$ have the same minimal
generators.  It follows that 
\[
I_\Delta = \bigcap_{i=1}^p \mathfrak{p}_{F_i} = \bigcap_{i=1}^p
\mathfrak{p}_{c,t}(m_i) \cdot R.
\]
Therefore, Theorem \ref{decompofh} shows that
\[
I_\Delta = I_{c-1,t}(\underline{g}) \cdot R,
\]
which means in particular that $I_\Delta$ defines a scheme in $\proj{n}$
which is a cone over $Z_{c-1,t}(\underline{g}) \subset \proj{c-2+2t}$. 
Using Lemma \ref{fax-9.1} again, we observe that $\mathfrak{p}$ is a minimal
prime of $I_{\partial \Delta}$ if and only if it is generated by $c$ of the
variables $x_0,\dots,x_n$ and contains exactly one minimal prime ideal of
$I_\Delta = I_{c-1,t}(\underline{g}) \cdot R$.  Hence Theorem
\ref{maingorthm} provides $I_{\partial \Delta} = J_c(\underline{h})$.  Since
the graded Betti numbers of $K[\partial \Delta] = R/J_c(\underline{h})$ were
computed in Theorem \ref{main-thm-paper}, the proof is complete.
\end{proof}

\begin{remark}
It is interesting to observe that the passage from $\Delta$ to $\partial
\Delta$ has an interpretation using Gorenstein liaison.  Indeed, the proof
above shows that $I_{\partial \Delta} = I_\Delta + I_Y$ where $Y$ is the
residual of $I_\Delta = I_{c-1,t}(\underline{g})$ in $I_{G_{c-1,s+1,t}}$.
\end{remark}

Using Corollary \ref{fax-8.14} instead of Theorem \ref{main-thm-paper}, the
proof above provides:

\begin{corollary}\label{fax-9.7}
Let $P$ be a simplicial polytope with $h$-vector $\underline{h}$ as given in
Corollary \ref{fax-8.14}.  Then the shape of the minimal free resolution  of
$K[\Delta(P)]$ is the one described in Corollary \ref{fax-8.14} if $K$ has
characteristic zero.
\end{corollary}

\begin{remark}
A simplicial $d$-polytope $P$ is called {\em stacked} if it admits a
triangulation $\Gamma$ which is a $(d-1)$-tree, i.e.\ $\Gamma$ is a
shellable
$(d-1)$-dimensional simplicial complex with $h$-vector $(1,c-1)$ (cf.\
\cite{HL}, Corollary 1.3).  The graded Betti numbers of stacked polytopes
have been computed by Hibi and Terai in \cite{H-T} (cf.\ also \cite{HL},
Theorem 3.3).  Since the $h$-vector of a stacked $d$-polytope with $s+c$
vertices is $(1,c,\dots,c,1)$, this result may be considered as a special
case of Corollary \ref{fax-9.7} (with $t = 1$) if $s \neq 3$. 
\end{remark}


\section{Final comments}

\begin{remark}\label{rmk-on-pf}
Our approach in Theorem \ref{maingorthm} (and hence Theorem
\ref{main-thm-paper}(b)) is similar to that of
\cite{GM5} in that we construct our linking scheme $X$ to be a generalized
stick figure, thus automatically guaranteeing that any choice of an
arithmetically Cohen-Macaulay subscheme $W \subset X$ will be geometrically
linked to its residual, so we can add the linked ideals to produce our
desired arithmetically Gorenstein scheme.  There are two important
differences between our approach and that of \cite{GM5}.  First, in
\cite{GM5} it was enough to choose $X$ to be a complete intersection (in
fact, $X$ had codimension two and so was forced to be a complete
intersection).  In our situation, complete intersections do not give all
possible Hilbert functions! (See Example \ref{ci-not-enough}.)  Also, we had
to apply the ``sums of linked ideals'' method more than once. Second, in
\cite{GM5}, as in most applications of liaison, the authors started with the
scheme $W$ and found  the ``right''
$X$ to geometrically link it and produce the desired residual.  In our case we
start with a very reducible $X$ and show that there must exist the ``right''
$W$ inside it!

We should also remark that our method is quite different from that of
Harima \cite{harima}.  In that paper the author heavily uses the Artinian
property.  He forms a geometric link of two finite sets of points using a
complete intersection, takes the sum of the linked ideals, and considers an
ideal quotient on the resulting Artinian ideal which modifies the Hilbert
function in the desired way.  To mimic that approach we would have to first
construct a larger Gorenstein scheme $Y'$ by linking with a complete
intersection, find a large number of components of $Y'$ lying on a hyperplane
and forming an arithmetically Cohen-Macaulay union of linear spaces with
certain Hilbert function, remove these components, repeat this process a
certain number of times, and show that the result of this procedure preserves
the Gorenstein property.  It is possible that such an approach would also
give a construction.  With the advantage of hindsight, though, something like
this does happen.  See Remark \ref{compare-harima-2}.
\end{remark}

\begin{example}\label{ci-not-enough} As mentioned above, if we try to apply
the approach of this paper but use only complete intersections for our links,
we cannot obtain all SI-sequences.  For example, we consider Remark 3.5 of
\cite{GKS}.  There, the authors say, they cannot use their method (basically
sums of certain linked ideals using complete intersection links) to obtain
Gorenstein rings with Hilbert function
\[
1 \ \ 4 \ \ 10 \ \ n \ \ 10 \ \ 4 \ \ 1
\]
where $14 \leq n \leq 20$, although they remark that it is known how to
construct them for $n=20$ using other methods.  Here we show that the cases
$n=20$ and $n=19$ cannot be constructed as the sum of complete intersection
linked ideals, but the other cases {\em can} be so constructed.  Of course all
of them can be constructed by our method with Gorenstein links.

Now we will use the notation $h_{G}(t)$ for the entry in degree $t$ of the
$h$-vector of $G$, and we 
will use the same notation as in Remark \ref{codim-c-1-to-codim-c} for the
linked schemes.  We know that we want a codimension 4 Gorenstein scheme $G_4$
with regularity 7, so the complete intersection $G_3$ would have to have
regularity 8.  Hence the sum of the degrees of the generators must be 10.
Furthermore, one can check that because we seek $\Delta h_{G_4} (2) = 6$, we
must have $h_{G_3}(2) = 6$.  That is, $G_3$ has no quadric generator.  Hence
$G_3$ is the complete intersection of two cubics and a quartic.
Furthermore, the linked scheme $Z_3$ must also have $h_{Z_3}(2) = 6$ and
$h_{Y_3}(2) = 6$.  Hence the regularity of $Z_3$ is at most 5.  We get the
following diagram:

\bigskip
\begin{center}

\begin{tabular}{r|ccccccccccccc}
&\multicolumn{3}{l}{degree:}  \\
 &$0$ &$1$ & $2$ & $3$ & $4$ &$5$ & $6$ & $7$ & $8$ \\ \hline
$G_3$  & $1$ & $3$ & $6$ & $8$ & $8$ & $6$ & $3$ & $1$ & $0$  \\
$Z_3$ & $1$ & $3$ & $6$ & $a$ & $b$ & $0$ \\
$Y_{3}$ & $ 1$ & $3$ & $6$ &
$8-b$ & $8-a$ & $0$ \\
$\Delta G_{4}$ &  $ 1$ & $ 3$ & $ 6$
&
$n-10$ & $ 10-n$ & $ -6$ & $
-3$ & $ -1$ & $0$ \\
$G_{4}$ &  $ 1$ & $ 4$ & $ 10$
& $ n$ & $ 10$ & $ 4$ &
$
1$ & $ 0$ & $0$
\end{tabular}
\end{center}

\medskip

\noindent where $a$ and $b$ still have to be determined.  As a result, we have
\[
n-10 = (8-b) +a-8 = a-b .
\]
Since clearly $a \leq 8$, we get that $n=20$ and $n=19$ cannot be achieved in
this way.  Furthermore, any $n$ in the range $14 \leq n \leq 18$ {\em can} be
constructed in this way.  It is enough to take $Z_3$ to be a set of $n$ general
points in $\proj{3}$ (hence $a= n-10$ and $b=0$) and $G_3$ the complete
intersection of two cubics and one quartic containing $Z_3$.
\end{example}

\begin{remark}\label{truncation}
In the process of thinking about the construction of the arithmetically
Gorenstein schemes in Theorem \ref{maingorthm}, we were led to
consider the following question (as a special case), which is of independent
interest.  Suppose that $X$ is an \acm union of lines in projective space.
Is it possible to remove lines from $X$ one by one such that at each step the
union of the remaining lines is arithmetically Cohen-Macaulay, and the Hilbert
function of the general hyperplane section is the truncated Hilbert function
in the sense of \cite{GKR}?  At first glance the answer would seem to be yes
(and we believe  it to be yes) since the hyperplane section is a finite set
of points, and it is known that one can take away a point at a time giving
truncated Hilbert functions.  But the problem is that removing the
corresponding line from $X$ does not guarantee that the remaining union of
lines will be arithmetically Cohen-Macaulay.

For example, let $X$ be the union of three lines $A$, $B$ and $C$ such that
$B$ meets each of $A$ and $C$ but $A$ and $C$ are disjoint.  This is \acm (of
arithmetic genus 0).  Clearly removing either $A$ or $C$ preserves the
property of being \acm while removing $B$ does not.
\end{remark}


\end{document}